\documentclass[12pt]{article}
\usepackage[utf8]{inputenc}
\usepackage[T1]{fontenc}
\usepackage[a4paper,margin=1in]{geometry}
\usepackage{amsmath,amssymb,amsthm,mathtools,bm,mathrsfs}
\usepackage{graphicx}
\usepackage{booktabs,array}
\usepackage{float}
\usepackage{enumitem}
\usepackage{hyperref}
\usepackage{array}
\usepackage{ragged2e}

\newtheorem{theorem}{Theorem}

\newtheorem{conjecture}{Conjecture}
\newtheorem{remark}{Remark}

\title{Young Measure Based Quantum Linear Programming Algorithms for Nonlinear/Stochastic Multiscale Partial Differential Equations and Homogenization }
\author{Siqi Chen, Shi Jin, Lei Zhang}
\date{\today}

\begin{document}
\raggedbottom
\setlength{\textfloatsep}{8pt plus 2pt minus 2pt}
\setlength{\floatsep}{6pt plus 2pt minus 2pt}
\setlength{\intextsep}{6pt plus 2pt minus 2pt}
\maketitle

\begin{abstract}
We study quantum algorithms for nonlinear and stochastic homogenization via a Young-measure based linear programming (LP) formulation, which lifts the nonlinear problem to a linear one in higher dimensions by treating the microscale, the gradient, and possible random variables as independent variables, thereby capturing effective macroscopic quantities without directly resolving fine-scale oscillations. The resulting LP is large but structured, and its high-dimensional nature creates regimes in which quantum LP solvers outperform direct classical solvers: in the deterministic setting, polynomial quantum speedup arises when moderate homogenized accuracy suffices; in the stochastic setting, encoding all random realizations simultaneously in a single LP yields a quantum square-root reduction in stochastic sampling cost that grows with the number of random variables. Regularity or sparsity of the Young measure may further extend these advantages to fine-scale accuracy. Numerical experiments on one- and two-dimensional benchmarks confirm the correctness of the Young-measure LP formulation.
\end{abstract}


\setlength{\textfloatsep}{8pt plus 2pt minus 2pt}
\setlength{\floatsep}{6pt plus 2pt minus 2pt}
\setlength{\intextsep}{6pt plus 2pt minus 2pt}

\section{Introduction}

Many physical, biological, and engineering systems are governed by partial differential equations (PDEs) whose coefficients vary on spatial or temporal scales much smaller than the scale at which the system is observed. Examples include flow through porous media, heat conduction in composite materials, elastic deformation of heterogeneous media, wave propagation in complex environments, and transport in random or disordered structures.  Homogenization theory provides a systematic framework for replacing the original heterogeneous problem by an effective macroscopic model that captures the large-scale behavior of the solution \cite{bensoussan1978asymptotic, jikov1994homogenization, cioranescu1999introduction}. They also guide the design of multiscale numerical methods for PDEs with highly oscillatory  coefficients \cite{hou1997multiscale, efendiev2009multiscale, weinan2003heterognous}. 

In the classical periodic setting, the microscopic medium is modeled by a repeating unit cell, allowing effective equations to be derived through asymptotic expansions, two-scale convergence, or variational methods. The homogenized coefficients describe the macroscopic response of the oscillatory medium. However, many realistic media are nonperiodic and exhibit disorder, randomness, intermittency, or long-range correlations. This motivates stochastic homogenization, where coefficients are modeled as stationary random fields rather than periodic functions \cite{kozlov1979averaging, papanicolaou1981boundary, gloria2016quantitative}, and studies when equations with random microscopic coefficients converge, at large scales, to deterministic effective equations. Under assumptions such as stationarity, ergodicity, coercivity, and suitable compactness or mixing conditions, solutions often converge to a deterministic homogenized limit. Effective coefficients are obtained not from periodic cell problems, but from corrector problems on probability spaces or large finite domains \cite{papanicolaou1981boundary, armstrong2016large, gloria2016quantitative}. Nonlinearity further complicates the theory, since effective behavior is not determined by simple coefficient averaging. Oscillations in both the coefficients and solution gradients may interact, and weak convergence is generally insufficient for passing to nonlinear limits. Thus, nonlinear homogenization must account for microscopic oscillations, loss of compactness, and nonlinear material response \cite{dalmaso1993introduction, braides2002gamma, evans1992periodic}.

Several tools support nonlinear and stochastic homogenization, including variational methods, $\Gamma$-convergence, compensated compactness, monotonicity methods, subadditive ergodic theorems, corrector estimates, concentration inequalities, and quantitative homogenization \cite{avellaneda1987homogenization, kozlov1979averaging, braides2002gamma, kenig2012convergence,armstrong2016large, gloria2016quantitative}. In this work, we emphasize Young measures as a central framework for describing nonlinear oscillations and randomness-induced limiting distributions that are not captured by weak convergence alone but are essential for passing nonlinear quantities to the limit \cite{young1969lectures, ball1989version, pedregal1997parametrized}. The goal is to connect variational structure, stochastic heterogeneity, Young measure descriptions, and numerical approximation in order to identify effective macroscopic laws without fully resolving every microscopic realization or scale. By retaining distributional information beyond averages, Young measure formulations provide a bridge between rigorous nonlinear and stochastic homogenization theory and computational approaches such as discretization, relaxation, linear programming, adaptive atomic approximations, and neural-network parameterizations \cite{carstensen2000numerical, karimi2025approximating}. A key computational challenge is that Young measures are probability-valued maps encoding distributions over gradients, phases, microscopic states, and random environments, so their approximation is intrinsically high-dimensional across physical, phase, and probability spaces.

Quantum computing for differential equations aims to reduce the cost of high-dimensional PDE computation by exploiting techniques such as quantum linear algebra and Hamiltonian simulation. For linear PDEs, standard approaches either discretize the equation and apply quantum linear-system solvers, or embed the resulting possibly nonunitary dynamics into Hamiltonian simulation frameworks, for example through Schrödingerization \cite{harrow2009quantum,berry2017quantum,jin2023quantum,jin2024schrodingerization}. For linear multiscale problems with scale separation, Hu et. al developed quantum algorithms for elliptic, parabolic, and hyperbolic PDEs based on either direct oscillatory models or homogenized equations, combined with higher-dimensional liftings and Schrödingerization \cite{hu2024quantum}. Balazi et. al proposed a quantum-enhanced numerical homogenization method based on the Localized Orthogonal Decomposition, achieving logarithmic scaling in fine-scale resolution via quantum local solvers \cite{balazi2026quantum}. Nonlinear PDEs are more challenging because quantum evolution is intrinsically linear. This motivates a Young-measure formulation \cite{jin2026young}, in which oscillatory, singular, or uncertain nonlinear solutions are represented by probability-valued measures. Under this representation, certain nonlinear PDE problems can be reformulated as linear programming (LP) problems with linear objectives and constraints. Since Young measures are defined on high-dimensional physical, state, and possibly random-parameter domains, their full computation suffers from the curse of dimensionality. Quantum linear-programming methods may therefore offer polynomial advantages for recovering the full Young measure, although such advantages may be less relevant when only the expected weak solution is required \cite{jin2026young,apers2023quantum,augustino2023quantum,kerenidis2020quantum}.

In this paper, we first introduce nonlinear heterogeneous and random homogenization models and formulate them using Young measures. We then show how the measure-valued formulation transforms the nonlinear PDE into linear constraints involving normalization, barycenter consistency, compatibility, and weak force balance. Next, we describe the resulting discrete Young-measure linear programs for deterministic periodic and random media. We then discuss quantum algorithms for solving these LPs, focusing on the quantum central path method and comparing its complexity with that of direct classical solvers and classical Young-measure LP approaches. 

The quantum advantage analysis is detailed in Sections~\ref{subsec:qcp-classical-comparison} and~\ref{subsec:qcp-regimes}. The baseline is a direct classical solver (DCS), which must resolve all $\varepsilon$-scale oscillations at cost $\widetilde O(\varepsilon^{-d})$ with $d$ being the spatial dimension. The Young-measure LP, by contrast, treats the microscale $y$, the gradient $\xi$, and the random environment $\omega$ as independent discretization variables, and its size grows with the target accuracy $\delta$ rather than with $\varepsilon$. Under first-order algebraic discretization, which requires minimal regularity assumptions, two rigorous regimes of quantum advantage over DCS are identified (see Table~\ref{tab:quantum-advantage}):
\begin{enumerate}
\item \emph{Deterministic nonlinear homogenization.} Recovering the homogenized solution $u_0$ to accuracy $\delta=\varepsilon^\alpha$, where $\alpha\in(0,1)$ reflects the homogenization convergence rate, yields a polynomial QCP speedup over DCS whenever $\alpha<2d/(3d+4)$.
\item \emph{Stochastic homogenization at fine-scale accuracy.} With $r$ random variables and $N_\omega^r$ total random discretization points treated as a free parameter (independent of $\varepsilon$), QCP achieves a square-root reduction in stochastic sampling cost whenever $N_\omega^r\gtrsim\varepsilon^{-(d+2)}$; the speedup grows with $N_\omega^r$.
\end{enumerate}
We further propose in Conjecture~\ref{conj:fine-solution} that if the Young-measure density is analytic in $x$ and $y$, spectral discretization in those variables extends quantum advantage to fine-scale accuracy $\delta=\varepsilon$ for all $d\ge2$ in the deterministic setting, with an analogous improvement in the stochastic case.

\begin{table}[H]
\centering
\small
\renewcommand{\arraystretch}{1.4}
\begin{tabular}{@{}lllll@{}}
\toprule
Setting & Target accuracy & DCS cost & QCP cost & Advantage condition \\
\midrule
Deterministic & $\delta=\varepsilon^\alpha$, $\alpha\in(0,1)$
  & $\widetilde O(\varepsilon^{-d})$
  & $\widetilde O(R_1\,\varepsilon^{-\alpha(3d+4)/2})$
  & $\alpha < \dfrac{2d}{3d+4}$ \\[6pt]
Stochastic & $\delta=\varepsilon$, $N_\omega^r$ free
  & $\widetilde O(N_\omega^r\,\varepsilon^{-d})$
  & $\widetilde O(R_1\,N_\omega^{r/2}\,\varepsilon^{-(3d/2+1)})$
  & $N_\omega^r \gtrsim \varepsilon^{-(d+2)}$ \\
\bottomrule
\end{tabular}
\caption{Quantum advantage regimes for the Young-measure LP using  first-order algebraic discretization. Here $\varepsilon$ is the microscale, $d$ is the spatial dimension, $R_1$ bounds the central-path length, and $N_\omega^r$ is the total number of random discretization points treated as a free parameter. DCS denotes direct classical solver. See Sections~\ref{subsec:qcp-classical-comparison} and~\ref{subsec:qcp-regimes} for details.}
\label{tab:quantum-advantage}
\end{table}

\emph{Organizatin}: The rest of the paper is organized as follows. Section 2 reviews the deterministic periodic and stochastic homogenization settings for nonlinear elliptic equations. Section 3 develops the Young-measure formulation and its finite-dimensional LP discretization for both periodic and random media. Section 4 discusses the quantum central path algorithm for the resulting Young-measure LPs, derives the associated query and gate complexity bounds, and compares them with direct classical multiscale solvers and classical LP methods. Section 5 presents numerical experiments for one- and two-dimensional deterministic and stochastic benchmark problems. Section 6 concludes with a summary of the proposed framework and outlines several directions for future work.

\emph{Notation}: $\widetilde{O}(f(n))$ denotes complexity scaling as $O(f(n))$ up to logarithmic factors.

\section{Preliminaries on elliptic homogenization}
\label{sec:homogenization}

We first recall the homogenization framework used throughout the paper. The
deterministic setting describes periodic nonlinear elliptic equations with
rapidly oscillating coefficients, while the random setting describes
media whose constitutive laws depend on finitely many probabilistic states.
In both cases, the goal is to replace the heterogeneous problem by an
effective deterministic equation whose coefficient, or constitutive law,
captures the macroscopic response of the medium.

\subsection{Nonlinear periodic homogenization}
\label{subsec:nonlinear_periodic_homogenization}

We first consider deterministic periodic homogenization for nonlinear
monotone elliptic equations. Let \(\Omega\subset\mathbb{R}^d\) be a bounded
domain and let \(Y=(0,1)^d\) denote the periodicity cell. For
\(\varepsilon>0\), we study
\begin{equation}
\label{eq:pde_eps}
-\nabla\cdot a\!\left(\frac{x}{\varepsilon},\nabla u_\varepsilon(x)\right)
=
f(x)
\quad\text{in }\Omega,
\qquad
u_\varepsilon=0
\quad\text{on }\partial\Omega,
\end{equation}
where \(a(y,\xi)\) is \(Y\)-periodic in \(y\). We assume that \(a\) satisfies
the standard monotone \(L^p\) conditions: for some \(p>1\) and constants
\(c_0,c_1>0\),
\begin{align}
    a(y,\xi)\cdot \xi
    &\ge c_0|\xi|^p-c_1,
    \label{eq:monotone_coercivity}
    \\
    |a(y,\xi)|
    &\le c_1(1+|\xi|^{p-1}),
    \label{eq:monotone_growth}
    \\
    \big(a(y,\xi)-a(y,\eta)\big)\cdot(\xi-\eta)
    &\ge 0.
    \label{eq:monotone_monotonicity}
\end{align}
Under these assumptions, \((u_\varepsilon)\) is uniformly bounded in
\(W_0^{1,p}(\Omega)\). Hence, along a subsequence,
\begin{equation}
    u_\varepsilon \rightharpoonup u_0
    \quad\text{weakly in }W_0^{1,p}(\Omega).
\end{equation}
The homogenized limit \(u_0\) solves
\begin{equation}
\label{eq:hom_pde}
-\nabla\cdot a_{\mathrm{hom}}(\nabla u_0)
=
f
\quad\text{in }\Omega,
\qquad
u_0=0
\quad\text{on }\partial\Omega.
\end{equation}

The effective constitutive law \(a_{\mathrm{hom}}\) is defined through a
periodic cell problem. For each \(\xi\in\mathbb{R}^d\), find
\(w_\xi\in W^{1,p}_{\#}(Y)\), with \(\int_Y w_\xi\,dy=0\), such that
\begin{equation}
\label{eq:cell_problem}
-\nabla_y\cdot
a\big(y,\xi+\nabla_y w_\xi(y)\big)
=
0
\quad\text{in }Y.
\end{equation}
Then
\begin{equation}
\label{eq:ahom_definition}
a_{\mathrm{hom}}(\xi)
=
\int_Y
a\big(y,\xi+\nabla_y w_\xi(y)\big)\,dy.
\end{equation}
Here \(W^{1,p}_{\#}(Y)\) denotes the Sobolev space of \(Y\)-periodic
functions with \(Y\)-periodic weak derivatives. The zero-average condition
fixes the additive constant in the cell corrector. This formulation is
standard in periodic homogenization and two-scale convergence; see, for
example, \cite{bensoussan1978asymptotic, jikov1994homogenization,
allaire1992homogenization, dal_maso1993introduction}.

\begin{remark}
In one spatial dimension, the cell formula reduces to explicit effective
laws. For the linear scalar case \(a(y,\xi)=\alpha(y)\xi\), one obtains the
harmonic mean
\begin{equation}
\label{eq:linear_1d_hom}
    \alpha_{\mathrm{hom}}
    =
    \left(\int_0^1 \alpha(y)^{-1}\,dy\right)^{-1}.
\end{equation}
For the nonlinear example \(a(y,\xi)=k(y)\xi^3\), the effective coefficient satisfies
\begin{equation}
\label{eq:nonlinear_1d_hom}
    k_{\mathrm{hom}}
    =
    \left(\int_0^1 k(y)^{-1/3}\,dy\right)^{-3}.
\end{equation}
\end{remark}

\subsection{Stochastic homogenization}
\label{subsec:random_homogenization}

We also consider homogenization in random media. Let
\((\Theta,\mathcal{F},\mathbb{P})\) be a probability space equipped with a
measure-preserving ergodic group action
\((\tau_x)_{x\in\mathbb{R}^d}\).
The constitutive law is described by a random monotone operator
\[
a:\Theta\times\mathbb{R}^d\to\mathbb{R}^d,
\qquad
(\omega,\xi)\mapsto a(\omega,\xi),
\]
which is assumed to be stationary in the sense that
\[
a(\tau_x\omega,\xi)
\]
represents the material response at the spatial position \(x\).

For \(\varepsilon>0\), we consider the stochastic boundary-value problem
\begin{equation}
\label{eq:sto_pde}
-\nabla\cdot
a\!\left(
\tau_{x/\varepsilon}\omega,
\nabla u_\varepsilon(x,\omega)
\right)
=
f(x)
\quad\text{in }\Omega,
\qquad
u_\varepsilon=0
\quad\text{on }\partial\Omega.
\end{equation}
Here the factor \(x/\varepsilon\) describes rapidly oscillating random
heterogeneities at the microscopic scale \(\varepsilon\).

Under the standard assumptions of stationarity, ergodicity,
uniform monotonicity, and polynomial growth,
the family \(u_\varepsilon(\cdot,\omega)\) converges,
for \(\mathbb{P}\)-almost every \(\omega\),
to a deterministic limit \(u_0\) solving the homogenized problem
\begin{equation}
\label{eq:sto_hom_limit}
-\nabla\cdot
a_{\mathrm{hom}}(\nabla u_0)
=
f
\quad\text{in }\Omega,
\qquad
u_0=0
\quad\text{on }\partial\Omega.
\end{equation}
The effective operator \(a_{\mathrm{hom}}\) is deterministic as a
consequence of ergodicity; see
\cite{kozlov1979averaging,papanicolaou1979boundary,
jikov1994homogenization,armstrong2016quantitative}.

In the finite-state setting used later in this paper, the random parameter
takes values
\[
\omega_i\in\Theta,
\qquad i=1,\ldots,N,
\]
with probabilities
\[
p_i>0,
\qquad
\sum_{i=1}^N p_i=1.
\]
Each state \(\omega_i\) corresponds to a constitutive law
\[
a_i(\xi)
:=
a(\omega_i,\xi).
\]
Accordingly, the random medium is represented by the finite ensemble
\begin{equation}
\label{eq:finite_random_medium}
\big\{(a_i,p_i)\big\}_{i=1}^N.
\end{equation}

To approximate the stochastic homogenized operator, we introduce in the
next section a Young-measure formulation together with a finite-dimensional
linear-programming discretization adapted to the finite-state random
setting.

\begin{remark}
Classical periodic homogenization is treated in
\cite{bensoussan1978asymptotic,jikov1994homogenization,
allaire1992homogenization}. Nonlinear monotone homogenization and
\(\Gamma\)-convergence methods are discussed in
\cite{dal_maso1993introduction,braides2002gamma}. Stochastic
homogenization originates in the works
\cite{kozlov1979averaging,papanicolaou1979boundary}
and has been further developed in modern quantitative forms in
\cite{gloriaotto2011,gloriaotto2012,gloriaotto2015,gloriaotto_corrector,armstrong2016quantitative}.
\end{remark}

\begin{remark}[One-dimensional nonlinear random example]
Consider
\[
W(\omega,\xi)
=
\frac1p c(\omega)|\xi|^p,
\qquad
a(\omega,\xi)
=
c(\omega)|\xi|^{p-2}\xi,
\]
with \(p>1\) and
\[
0<c_0\le c(\omega)\le c_1.
\]
The corresponding one-dimensional equation is
\[
-\frac{d}{dx}
\left(
c(\tau_{x/\varepsilon}\omega)
|u_x^\varepsilon|^{p-2}u_x^\varepsilon
\right)
=
f(x).
\]

In the zero-forcing case, the flux is constant and the homogenized constitutive law takes the form
\[
J
=
c_{\mathrm{hom}}
|q|^{p-2}q,
\]
where the effective coefficient is
\[
c_{\mathrm{hom}}
=
\left(
\mathbb E\!\left[
c(\omega)^{-\frac1{p-1}}
\right]
\right)^{-(p-1)}.
\]

For a two-phase random medium,
\[
c(\omega)=
\begin{cases}
c_1,
& \text{with probability }\theta,
\\
c_2,
& \text{with probability }1-\theta,
\end{cases}
\]
one obtains
\[
c_{\mathrm{hom}}
=
\left(
\theta c_1^{-\frac1{p-1}}
+
(1-\theta)c_2^{-\frac1{p-1}}
\right)^{-(p-1)}.
\]

When \(p=2\), this reduces to the classical harmonic mean formula.
\end{remark}

\section{Young-measure formulation and LP implementation}

This section presents a measure-valued formulation for nonlinear homogenization problems in both periodic and stochastic settings, together with the finite-dimensional linear programming problems used in the numerical implementation. The central idea is to replace rapidly oscillatory gradients by parametrized Young measures that encode the local distribution of microscopic states~\cite{young1969lectures, ball1989version, pedregal1997parametrized}. In both the periodic and stochastic frameworks, the resulting constraints: normalization, barycenter consistency, weak equilibrium, and admissibility conditions, are {\it linear} in the measures. Consequently, after discretization of the spatial variables and gradient state space, the homogenization problem becomes a large but structured linear program. The formulation provides a unified framework for deterministic periodic media and stochastic heterogeneous materials while avoiding direct resolution of fine microscopic scales.
\subsection{Nonlinear periodic homogenization}

Let $\Omega\subset\mathbb R^d$ be bounded and let $Y=(0,1)^d$ denote the periodic reference cell. Consider the nonlinear elliptic problem
\begin{equation}
-\nabla\cdot a\!\left(\frac{x}{\varepsilon},\nabla u^\varepsilon(x)\right)=f(x)
\qquad \text{in }\Omega,
\end{equation}
subject to the boundary condition
\[
u^\varepsilon=g
\qquad \text{on }\partial\Omega.
\]
Here
\[
a:Y\times\mathbb R^d\to\mathbb R^d
\]
is assumed to be $Y$-periodic in the first variable.

To describe microscopic oscillations in the gradients, we introduce a parametrized Young measure
\begin{equation}
\mu_{x,y}\in\mathcal P(\mathbb R^d),
\qquad (x,y)\in\Omega\times Y,
\end{equation}
which represents the local probability distribution of microscopic gradients at macroscopic position $x$ and microscopic phase $y$. The normalization condition is imposed pointwise:
\begin{equation}
\int_{\mathbb R^d} d\mu_{x,y}(\xi)=1
\qquad \text{for a.e. }(x,y)\in\Omega\times Y.
\end{equation}

The admissible unknowns are
\[
u\in W^{1,p}(\Omega),
\qquad
\mu_{x,y}\in\mathcal P(\mathbb R^d),
\]
subject to the following constraints:
\begin{align}
&\mu_{x,y}\ge0,
\\
&\int_{\mathbb R^d} d\mu_{x,y}(\xi)=1
\qquad \text{for a.e. }(x,y),
\\
&\nabla u(x)
=
\int_Y\int_{\mathbb R^d}
\xi\, d\mu_{x,y}(\xi)\,dy
\qquad \text{for a.e. }x,
\\
&-\nabla\cdot
\left(
\int_Y \int_{\mathbb R^d}
a(y,\xi)\,
d\mu_{x,y}(\xi)\,dy
\right)
=
f(x)
\qquad \text{in } \Omega .
\end{align}


\begin{remark} 
When the constitutive law derives from an energy density,
\[
a(y,\xi)=\partial_\xi W(y,\xi),
\]
the problem admits a variational interpretation, the associated measure-valued functional is
\begin{equation}
\mathcal J_{\mathrm{per}}(u,\mu)
=
\int_\Omega
\int_Y
\int_{\mathbb R^d}
W(y,\xi)\,
d\mu_{x,y}(\xi)\,dy\,dx
-
\int_\Omega f(x)u(x)\,dx.
\end{equation}
\end{remark}

Importantly, all terms involving $\mu$ are {\it linear} in the measure. After discretization, the formulation therefore reduces to a linear programming problem  in the discrete measure masses.

For presentation purposes, it is convenient to introduce the marginal measure
\begin{equation}
\widehat\mu_x(\xi)
:=
\int_Y \mu_{x,y}(\xi)\,dy,
\end{equation}
which represents the effective distribution of microscopic gradients at the macroscopic point $x$.

\subsection{Stochastic homogenization}

We now consider the nonlinear random elliptic problem
\begin{equation}
-\nabla\cdot
a\!\left(\tau_{x/\varepsilon}\omega,\nabla u^\varepsilon(x)\right)
=
f(x)
\qquad \text{in }\Omega,
\end{equation}
with boundary condition
\[
u^\varepsilon=g
\qquad \text{on }\partial\Omega.
\]

Here $(\Theta,\mathcal F,\mathbb P)$ is a probability space of environments,
\[
\tau_z:\Theta\to\Theta
\]
is a measure-preserving ergodic group action, and
\[
a:\Theta\times\mathbb R^d\to\mathbb R^d
\]
is the nonlinear random constitutive law. 

For almost every macroscopic point $x\in\Omega$, we introduce a probability measure
\begin{equation}
\mu_x\in\mathcal P(\Theta\times\mathbb R^d),
\end{equation}
interpreted as the local joint distribution of the random environment and microscopic gradient state. The admissibility conditions become
\begin{align}
&\mu_x\ge0,
\qquad 
\int_{\mathbb{R}^d} d\mu_x(\omega,\xi) = d\mathbb{P}(\omega) \quad \mathrm{for~a.e.}~x,
\\
&\nabla u(x)
=
\int_{\Theta\times\mathbb R^d}
\xi\,d\mu_x(\omega,\xi)
\qquad \text{for a.e. }x,
\\
&-\nabla\cdot
\left(
\int_{\Theta\times\mathbb R^d}
a(\omega,\xi)\,
d\mu_x(\omega,\xi)
\right)
=
f
\qquad \text{in } \Omega,
\\
&u=g
\qquad \text{on }\partial\Omega.
\end{align}

\begin{remark}
When $a(\omega,\xi)=\partial_\xi W(\omega,\xi)$, the problem again admits a variational formulation, the corresponding measure-valued objective functional is
\begin{equation}
\mathcal J_{\mathrm{rand}}(u,\mu)
=
\int_\Omega
\int_{\Theta\times\mathbb R^d}
W(\omega,\xi)\,
d\mu_x(\omega,\xi)\,dx
-
\int_\Omega f(x)u(x)\,dx.
\end{equation}
\end{remark}

As in the periodic case, the objective and all admissibility constraints are linear in the measures. For numerical implementation, the probability space is approximated by finitely many environments
\[
\omega_i,
\qquad
i=1,\dots,N_\omega,
\]
with associated probabilities
\[
p_i>0,
\qquad
\sum_{i=1}^{N_\omega}p_i=1.
\]
In contrast to the periodic setting, there is no explicit fast variable $y$. The unknown measure family is therefore indexed only by the macroscopic point and the random state:
\[
\mu_{x,i}\ge0.
\]

The constraints become
\begin{align}
&\int_{\mathbb R^d}
d\mu_{x,i}(\xi)=p_i
\qquad
\text{for a.e. }x,
\quad
i=1,\dots,N_\omega,
\\
&\nabla u(x)
=
\sum_{i=1}^{N_\omega}
\int_{\mathbb R^d}
\xi\,d\mu_{x,i}(\xi)
\qquad
\text{for a.e. }x,
\\
&-\nabla\cdot
\left(
\sum_{i=1}^{N_\omega}
\int_{\mathbb R^d}
a_i(\xi)\,
d\mu_{x,i}(\xi)
\right)
=
f
\qquad \text{in } \Omega,
\end{align}
where
\[
a_i(\xi)=a(\omega_i,\xi).
\]

Equivalently, one may introduce normalized probability measures
\[
\nu_{x,i}:=\frac{\mu_{x,i}}{p_i}
\in\mathcal P(\mathbb R^d),
\]
so that
\[
\mu_{x,i}=p_i\nu_{x,i}.
\]

The global marginal distribution used for visualization is
\begin{equation}
\widehat\mu_x(\xi)
=
\sum_{i=1}^{N_\omega}\mu_{x,i}(\xi)
=
\sum_{i=1}^{N_\omega}
p_i\nu_{x,i}(\xi).
\end{equation}

The discrete random objective functional becomes
\begin{equation}
\mathcal J_{\mathrm{rand}}(u,\mu)
=
\int_\Omega
\sum_{i=1}^{N_\omega}
\int_{\mathbb R^d}
W_i(\xi)\,
d\mu_{x,i}(\xi)\,dx
-
\int_\Omega f(x)u(x)\,dx,
\end{equation}
where
\[
W_i(\xi)=W(\omega_i,\xi).
\]

\begin{remark}
\label{rem:stochastic-two-scale-ym}
The formulation here is a \emph{mean} (annealed) formulation corresponding to stochastic two-scale convergence in the mean \cite{bourgeat1994stochastic}, in which the test functions $\varphi(\tau_{x/\varepsilon}\omega,x)$ are averaged over $P$; the \emph{quenched} variant in \cite{zhikov2006homogenization} instead fixes a realization $\omega_0$ and evaluates at $\tau_{x/\varepsilon}\omega_0$, with both notions coinciding under ergodicity. Heida et.al \cite{heida2021stochastic} unify these via \emph{stochastic two-scale Young measures} $\boldsymbol{\nu}=\{\nu_\omega\}_{\omega\in\Omega}$, which simultaneously encode the quenched pathwise behavior and the mean limit, with the Dirac case $\nu_\omega=\delta_{v(\omega)}$ recovering quenched from mean convergence. For the LP framework, the mean formulation encodes all $N_\omega^r$ environments in a single system of size $n_{\rm sto}\sim N_x^d N_y^d N_\omega^r N_\xi^d$, enabling the quantum square-root advantage in the stochastic dimension (Section~\ref{subsec:qcp-regimes}), whereas a quenched LP would solve $N_\omega^r$ independent problems at the same total cost but without this advantage.
\end{remark}

\subsection{Discrete LP implementation}

After discretization of the macroscopic variable and, when required, the microscopic phase, random state, and gradient variable, all unknowns are assembled into a vector
\begin{equation}
z=[u;F],
\qquad
F\ge0,
\end{equation}
where $F$ contains the discrete masses approximating the Young measures $\mu_{x,y}$ in the periodic setting or $\mu_{x,i}$ in the stochastic setting.

The resulting finite-dimensional optimization problem takes the form
\begin{equation}
\min_z c^\top z
\qquad
\text{subject to}
\qquad
A_{\mathrm{eq}}z=b_{\mathrm{eq}},
\qquad
F\ge0.
\end{equation}

The equality constraints enforce:
\begin{align}
&\text{local normalization or state-wise mass conservation},\\
&\text{consistency between barycenters and macroscopic gradients},\\
&\text{weak equilibrium or force balance},\\
&\text{boundary conditions on }u.
\end{align}


\begin{remark}
If the constitutive law admits a potential representation \(a=\partial_\xi W\), the vector \(c\) arises from the discretized microscopic energy and the LP corresponds to an energy minimization problem. Otherwise, the formulation remains linear in the measure variables, but the LP is interpreted primarily as a feasibility problem for admissible measure-valued equilibria.
\end{remark}

\begin{remark}
The Young-measure formulation does not eliminate computational complexity. Rather, it replaces direct fine-scale resolution of multiscale PDEs by a larger optimization problem over measure-valued variables indexed by macroscopic position, microscopic phase, random environment, and discretized gradient state.
\end{remark}

\begin{remark}
In the stochastic setting, additional admissibility conditions may be required to ensure that the measures arise from stationary random gradient fields. Such constraints may include stationarity, curl-free compatibility, moment bounds, or corrector-type conditions.
\end{remark}

\section{Quantum Algorithm for Young Measures}
\label{sec:quantum-algorithm-young-measures}

This section examines the application of the quantum central path (QCP) algorithm to the Young-measure formulation of the homogenization problem. After reformulating the discretized Young-measure model as a finite-dimensional linear program, we review the homogeneous self-dual embedding underlying QCP and derive the corresponding query and gate complexity bounds. These complexities are then compared with those of direct fine-scale PDE solvers and classical interior-point methods. The analysis identifies parameter regimes—including small $\varepsilon$, high-dimensional measure discretizations, and sparse LP structure—in which QCP may provide computational advantages, particularly when the Young-measure discretization is independent of the microscopic scale and the resulting LP remains sparse and well-conditioned.

\subsection{Complexity of the Quantum Central Path Algorithm}
\label{subsec:qcp-complexity}

We consider the finite-dimensional Young-measure linear program obtained from the discretization introduced in the previous section. Let
\[
z=[u;F],
\qquad
F\ge0,
\]
where \(u\) denotes the discrete macroscopic variable and \(F\) collects the nonnegative discrete Young-measure masses. The discrete problem takes the form
\begin{equation}
\label{eq:qcp-young-measure-lp}
    \min_{z}
    \quad
    c^\top z
    \qquad
    \text{subject to}
    \qquad
    Az=b,
    \qquad
    z\ge0.
\end{equation}
The matrix $A\in\mathbb R^{m\times n}$
collects the discrete normalization conditions, barycenter consistency relations, weak equilibrium equations, and boundary constraints.

Following the homogeneous self-dual embedding, the LP is reformulated in variables $w\in\mathbb R^d,
d=m+n+2$. 
The embedding is written as
\begin{equation}
\label{eq:qcp-self-dual-embedding}
    s(w)=Mw+q,
    \qquad
    w\ge0,
    \qquad
    s(w)\ge0,
\end{equation}
where \(M^\top=-M\). For a barrier parameter \(\mu>0\), the central path is characterized by
\begin{equation}
\label{eq:qcp-central-path}
    w_i s_i(w)=\mu,
    \qquad
    i=1,\dots,d.
\end{equation}

Equivalently, define the central-path potential
\begin{equation}
\label{eq:qcp-potential}
    f_\mu(w)
    :=
    \frac12
    \left\|
        w\circ s(w)-\mu\mathbf1
    \right\|_2^2.
\end{equation}
The central-path solution \(w(\mu)\) is the global minimizer of \(f_\mu\), satisfying
\begin{equation}
\label{eq:qcp-potential-minimum}
    f_\mu(w(\mu))=0.
\end{equation}

\begin{theorem}[Query and gate complexity of QCP]
\label{thm:qcp-complexity}
Assume that the homogeneous self-dual embedding is regular and that the central-path solution satisfies $|z^\star|_1 \le R_1$. Then the quantum central path algorithm computes an $\varepsilon$-optimal solution using $\widetilde O\left(\sqrt{m+n},\frac{R_1}{\varepsilon}\right)$ queries to an oracle evaluating $f_\mu$. In the standard gate model, without QRAM, the corresponding gate complexity is $\widetilde O\left(\sqrt{m+n},\mathrm{nnz}(A)\frac{R_1}{\varepsilon}\right)$, where $\mathrm{nnz}(A)$ denotes the number of nonzero entries of the LP constraint matrix.
\end{theorem}

Here we give a heuristic explanation for the theorem. 
The QCP algorithm associates each value of the barrier parameter \(\mu\) with a Schrödinger Hamiltonian
\begin{equation}
\label{eq:qcp-hamiltonian}
    H_\mu
    =
    -\frac{h}{2}\Delta + f_\mu(z),
\end{equation}
where \(h>0\) is a semiclassical parameter. Since \(f_\mu\) is minimized at the central-path point \(z(\mu)\), the ground state of \(H_\mu\) is concentrated near \(z(\mu)\).

Near \(z(\mu)\), we have the Taylor expansion
\begin{equation}
\label{eq:qcp-taylor}
    f_\mu(z)
    =
    \frac12
    (z-z(\mu))^\top
    \nabla^2 f_\mu(z(\mu))
    (z-z(\mu))
    +
    O(\|z-z(\mu)\|^3),
\end{equation}
where we have used \(f_\mu(z(\mu))=0\) and
\(\nabla f_\mu(z(\mu))=0\). Thus, locally, \(H_\mu\) is approximated by a multidimensional quantum harmonic oscillator. Its ground state is approximately Gaussian and localized around \(z(\mu)\).

To recover an \(\varepsilon\)-accurate classical point in dimension
\(d=m+n+2\), the wavepacket must be localized at scale
\(O(\varepsilon/\sqrt d)\) in each coordinate. Indeed, if the coordinatewise error is \(O(\eta)\), then the Euclidean error is \(O(\sqrt d\,\eta)\). Hence one needs
$\displaystyle \eta = O\!\left(\frac{\varepsilon}{\sqrt d}\right)$. The total length of the central path is controlled by the assumed bound $\|z^\star\|_1 \le R_1$. Therefore, the effective number of localization-scale steps required to transport the quantum state along the path is $\displaystyle
    \frac{R_1}{\varepsilon/\sqrt d}
    =
    \sqrt d\,\frac{R_1}{\varepsilon}$. Up to logarithmic factors hidden in \(\widetilde{O}(\cdot)\), this yields the query complexity
$\displaystyle \widetilde{O}\!\left(
        \sqrt d\,\frac{R_1}{\varepsilon}
    \right)$.

Substituting \(d=m+n+2\), and absorbing constants into the asymptotic notation, gives
\[
    \widetilde{O}\!\left(
        \sqrt{m+n}\,\frac{R_1}{\varepsilon}
    \right).
\]

Finally, in the standard gate model, each query to \(f_\mu\) requires evaluating
$ s(z)=Mz+q$. Since \(M\) is constructed from the original constraint matrix \(A\), this evaluation costs
$\widetilde{O}(\mathrm{nnz}(A))$
elementary gates. Multiplying this cost by the query complexity gives the overall gate complexity
\[
    \widetilde{O}\!\left(
        \sqrt{m+n}\,\mathrm{nnz}(A)\frac{R_1}{\varepsilon}
    \right).
\]

\subsection{Setup for the comparison with classical solvers}
\label{subsec:qcp-classical-comparison}

We compare the QCP complexity with direct classical solvers (DCS) for the multiscale PDE, using the Young-measure LP degrees of freedom and discretization error as the bridge between the two approaches.

\paragraph{Degrees of freedom.}
We consider both deterministic and stochastic Young-measure formulations. In the deterministic
periodic setting, let $x\in\Omega\subset\mathbb{R}^d$, $y\in Y=[0,1]^d$, and $\xi\in\mathbb{R}^d$.
We discretize $x$, $y$, and $\xi$ using $N_x^d$, $N_y^d$, and $N_\xi^d$ degrees of freedom
respectively. The discrete Young-measure variable is $\mu_{i,j,\ell}\ge0$, where
$i=1,\dots,N_x^d$, $j=1,\dots,N_y^d$, and $\ell=1,\dots,N_\xi^d$, giving
\[
n_{\rm det} \sim N_x^d N_y^d N_\xi^d.
\]

In the stochastic setting, we additionally let $\omega\in\Gamma\subset\mathbb{R}^r$ be
a random parameter with $r$ random variables, discretized on a tensor grid of $N_\omega^r$
points. The discrete Young-measure variable is $\mu_{i,j,k,\ell}\ge0$, where
$i=1,\dots,N_x^d$, $j=1,\dots,N_y^d$, $k=1,\dots,N_\omega^r$, and $\ell=1,\dots,N_\xi^d$,
giving
\[
n_{\rm sto} \sim N_x^d N_y^d N_\omega^r N_\xi^d.
\]

In both cases the LP is written abstractly as $\min_{z\in\mathbb{R}^n} c^\top z$ subject
to $A_{\rm LP}z=b$ and $z\ge0$, with $m\sim n$ for local Young-measure formulations.
The QCP query complexity is $\widetilde O(\sqrt{n}\,R_1/\delta)$.

\paragraph{Discretization error.}
We use first-order ($p=1$) discretization in $x$, $y$, and $\xi$, which requires minimal regularity assumptions on the Young-measure density:
\[
\mathrm{err}_{\rm LP}
\;\sim\;
N_x^{-1} + N_y^{-1} + N_\xi^{-1} + \varepsilon_{\rm alg}.
\]
Balancing at target accuracy $\delta$ gives $N_x\sim N_y\sim N_\xi\sim\delta^{-1}$ and
$\varepsilon_{\rm alg}\sim\delta$, so
\[
n_{\rm det}\sim\delta^{-3d},\qquad n_{\rm sto}\sim\delta^{-3d} N_\omega^r.
\]
The case of exponential (spectral) convergence in $x$ and $y$, which requires analytic
regularity of the Young-measure density in those variables
(see Section~\ref{subsec:qcp-regimes}), is deferred to Conjecture~\ref{conj:fine-solution}.

\paragraph{Desired accuracy.}
We distinguish two target accuracy regimes.
\begin{itemize}
\item \emph{Homogenized-solution regime.} The desired output is the homogenized solution
$u_0$. The natural accuracy threshold is
\[
\delta_{\rm hom} = \varepsilon^\alpha, \quad 0<\alpha<1,
\]
since the homogenization error satisfies $\|u_\varepsilon-u_0\|\sim\varepsilon^\alpha$.
This threshold is coarser than resolving the full oscillatory field.

\item \emph{Fine-solution regime.} The desired output is the full fine-scale solution
$u_\varepsilon$, with threshold $\delta_{\rm fine}=\varepsilon$.
Under $p=1$ algebraic discretization this regime does not yield quantum advantage
(see Conjecture~\ref{conj:fine-solution} for the spectral case).
\end{itemize}

\begin{remark}
\label{rem:norm-exponent-choice}
The discretization error above is a schematic complexity-level guide. The exponent $\alpha$ and the relevant norm depend on the model: for nonlinear monotone problems through $W^{1,p}$ or $L^p$ norms, and for stochastic homogenization on stationarity, ergodicity, and corrector estimates; see \cite{dal_maso1993introduction,braides2002gamma,evans1992periodic,kozlov1979averaging,armstrong2016quantitative,gloria2016quantitative}.
\end{remark}

\subsection{Possible quantum advantage regimes}
\label{subsec:qcp-regimes}

We now describe three regimes in which QCP may outperform DCS, using the setup from Section~\ref{subsec:qcp-classical-comparison}. Since $m\sim n$ for the sparse local Young-measure LPs considered here, the QCP query complexity reduces to $\widetilde O\!\bigl(\sqrt{n}\,R_1/\delta\bigr)$.

\paragraph{1. Nonlinear homogenization with homogenized accuracy ($\delta=\varepsilon^\alpha$, $0<\alpha<1$, $p=1$).}
We seek the homogenized solution $u_0$ to accuracy $\delta=\varepsilon^\alpha$ with $0<\alpha<1$.
This is the natural target in nonlinear homogenization, where the
homogenization error satisfies $\|u_\varepsilon-u_0\|\sim\varepsilon^\alpha$
with $\alpha<1$ in general norms (for instance $\alpha=1/2$ in $H^1$ for smooth linear
periodic problems, and possibly smaller for nonlinear settings). A direct
fine-scale solver must resolve all $\varepsilon$-scale oscillations regardless of the desired
accuracy for $u_0$, incurring cost
\[
\mathrm{DCS}_{\rm det}\sim\widetilde O(\varepsilon^{-d}).
\]

For the Young-measure LP with first-order ($p=1$) algebraic convergence in $x$, $y$, and $\xi$,
balancing all discretization errors at $\delta=\varepsilon^\alpha$ gives
$N_x\sim N_y\sim N_\xi\sim\varepsilon^{-\alpha}$, so $n_{\rm det}\sim\varepsilon^{-3\alpha d}$, and
\[
\mathrm{QCP}_{\rm det}
\sim
\widetilde O\!\left(
R_1\,\varepsilon^{-\alpha(2+3d/2)}
\right)
=
\widetilde O\!\left(
R_1\,\varepsilon^{-\alpha(3d+4)/2}
\right).
\]
QCP beats DCS when $\alpha(3d+4)/2 < d$, i.e.,
\[
\boxed{\alpha < \frac{2d}{3d+4}.}
\]
For $d=2$ this requires $\alpha<4/10=2/5$; for $d=3$, $\alpha<6/13$; for $d>>1$, $\alpha\lesssim 2/3$. Since $\alpha<1$ is the
natural regime for nonlinear homogenization, this condition can be satisfied
for sufficiently small $\alpha$ in all dimensions $d\ge1$.

\paragraph{2. Stochastic homogenization with fine-scale accuracy ($\delta=\varepsilon$, $p=1$).}
The second advantage regime arises in stochastic homogenization with $r$ random variables
$\omega\in\mathbb{R}^r$, where we seek the fine-scale solution at accuracy $\delta=\varepsilon$.
We treat the total number of random discretization points $N_\omega^r$ as a free parameter,
independent of the spatial accuracy $\varepsilon$, since the convergence rate in $\omega$ depends
on the smoothness of the solution in the random parameter and need not match the spatial rate.
A direct solver must handle each of the $N_\omega^r$ realizations at fine-scale cost $\varepsilon^{-d}$:
\[
\mathrm{DCS}_{\rm sto}\sim\widetilde O(N_\omega^r\,\varepsilon^{-d}).
\]
The stochastic Young-measure LP encodes all random realizations simultaneously in a single LP.
With $p=1$ algebraic convergence in $x$, $y$, $\xi$ and $\delta=\varepsilon$,
$N_x\sim N_y\sim N_\xi\sim\varepsilon^{-1}$, so
$n_{\rm sto}\sim\varepsilon^{-3d}\,N_\omega^r$, and
\[
\mathrm{QCP}_{\rm sto}
\sim
\widetilde O\!\left(
R_1\,\varepsilon^{-(3d/2+1)}\,N_\omega^{r/2}
\right).
\]
QCP beats DCS when $\varepsilon^{-(3d/2+1)}\,N_\omega^{r/2}\lesssim N_\omega^r\,\varepsilon^{-d}$, i.e.,
\[
\boxed{N_\omega^r \gtrsim \varepsilon^{-(d+2)}.}
\]
This condition depends only on the total stochastic dimension $N_\omega^r$ and not on $r$ or
$N_\omega$ individually: whenever the random environment is discretized finely enough that
$N_\omega^r\gg\varepsilon^{-(d+2)}$, quantum advantage holds. The quantum speedup factor
is $N_\omega^{r/2}$ — a square root in the stochastic dimension — which grows without bound
as $N_\omega^r$ increases, independently of $\varepsilon$. When the advantage condition is met,
the DCS cost is at least $\varepsilon^{-(2d+2)}$, while the QCP cost is
$\widetilde O(R_1\,\varepsilon^{-(3d/2+1)}\,N_\omega^{r/2})$, giving a genuine polynomial separation.

\paragraph{3. Conjecture: quantum advantage under spectral discretization.}
Above discussions use first-order algebraic convergence in $x$, $y$, and $\xi$. As noted in Remark~\ref{rem:higher-order}, one may improve convergence in $x$ and $y$ while keeping first-order in $\xi$, since regularity of $\nu_{x,y}$ in $\xi$ is not inherited from smoothness of $a$. In the extreme case where $\nu_{x,y}(\xi)$ is analytic in $x$ and $y$, spectral discretization gives $N_x\sim N_y\sim\log(1/\varepsilon)$ while $N_\xi\sim\varepsilon^{-1}$, reducing the LP size from $n\sim\varepsilon^{-3d}$ (all algebraic) to $n_{\rm det}\sim\varepsilon^{-d}|\log\varepsilon|^{2d}$ at fine-scale accuracy $\delta=\varepsilon$.

\begin{conjecture}[Quantum advantage under spectral discretization in $x,y$]
\label{conj:fine-solution}
Suppose $\nu_{x,y}(\xi)$ is analytic in $x$ and $y$, $\xi$ is discretized algebraically at $N_\xi\sim\varepsilon^{-1}$, and $R_1$ is polynomially bounded in $|\log\varepsilon|$. At fine-scale accuracy $\delta=\varepsilon$,
\[
n_{\rm det}\sim\varepsilon^{-d}|\log\varepsilon|^{2d},\qquad
\mathrm{QCP}_{\rm det}\sim\widetilde O\!\bigl(R_1\,\varepsilon^{-(1+d/2)}|\log\varepsilon|^{d}\bigr),
\]
which improves over DCS cost $\widetilde O(\varepsilon^{-d})$ for all $d\ge2$. For the stochastic case with $N_\omega^r$ free, the advantage condition becomes $N_\omega^r\gtrsim\varepsilon^{-(d-2)}$, trivially satisfied for all $d\ge2$ and $r\ge1$.
\end{conjecture}

\begin{remark}[Higher-order discretization]
\label{rem:higher-order}
If $\nu_{x,y}(\xi)$ is smooth in $x$ and $y$, one may use $p$-th order schemes in $x,y$ while keeping first-order in $\xi$, giving $N_x\sim N_y\sim\delta^{-1/p}$ and $N_\xi\sim\delta^{-1}$, so $n\sim\delta^{-(2d/p+d)}$ and QCP cost $\widetilde O(R_1\,\delta^{-(d/p+d/2+1)})$. The quantum advantage condition for the deterministic homogenization regime becomes $\alpha < 2dp/(2d+dp+2p)$, which approaches $2d/(d+2)$ as $p\to\infty$ --- matching the spectral limit of Conjecture~\ref{conj:fine-solution}. First-order discretization in $\xi$ is the natural default because regularity of $\nu_{x,y}$ in $\xi$ is not inherited from smoothness of $a(y,\xi)$: in the convex case the Young measure is a Dirac mass, while in the nonconvex case its support is determined by the quasiconvex hull geometry, which can be non-smooth even for analytic $a$. Higher-order convergence in $\xi$ requires additional structural information such as finite support or laminate structure.
\end{remark}

%

\section{Numerical Experiments}
We validate the proposed Young-measure LP formulation on several benchmark homogenization problems, including 1D deterministic periodic media, 1D stochastic media, and 2D anisotropic periodic media. The primary purpose of these experiments is to verify that the LP formulation correctly reproduces the expected homogenized solutions in linear, nonlinear and stochastic settings.

The numerical error is measured using the relative \(L^\infty\)-error
\begin{equation}
\label{eq:rel_err_def}
\frac{
\|u_{\mathrm{num}}-u_{\mathrm{exact}}\|_\infty
}{
\|u_{\mathrm{exact}}\|_\infty
}.
\end{equation}

Each benchmark uses a problem-dependent discretization consisting of:
\begin{itemize}
\item macroscopic spatial grids,
\item microscopic periodic or random-state discretizations,
\item gradient-state discretizations for the Young measures,
\item and the associated constitutive law.
\end{itemize}
The corresponding parameters are summarized in the tables preceding each experiment.

The experiments compare the LP-computed homogenized solution with analytically known homogenized solutions. In all cases considered below, the numerical results exhibit good agreement with the exact homogenized profiles and small relative errors, supporting the correctness and consistency of the Young-measure LP formulation.


\subsection{1D Deterministic Periodic Media}
\subsubsection{A linear case}
\begin{table}[H]
\centering
\small
\begin{tabular}{@{}ll@{}}
\toprule
Parameter & Value \\
\midrule
Macro grid & $50$ \\
Micro grid & $30$ \\
State grid $\xi$ & $201$, range $[-1,1]$ \\
Constitutive law & $W=\frac{1}{2}k(y)\xi^2$, $a=k(y)\xi$ (linear) \\
\bottomrule
\end{tabular}
\caption{A 1D deterministic linear: LP discretization setup.}
\label{tab:disc_1d_det_lin}
\end{table}

The oscillatory problem is $\displaystyle -\frac{d}{dx}\left(k\!\left(\frac{x}{\varepsilon}\right)u_\varepsilon'(x)\right)=1,
 u_\varepsilon(0)=u_\varepsilon(1)=0$, with periodic coefficient
$k(y)=2+\sin(2\pi y)$.

The effective coefficient is
$\displaystyle k_{\mathrm{hom}}=\left(\int_0^1\frac{1}{k(y)}\,dy\right)^{-1}=\sqrt{3}$,
and the homogenized problem
$\displaystyle -\frac{d}{dx}\left(k_{\mathrm{hom}}u_0'(x)\right)=1,
 u_0(0)=u_0(1)=0,
u_{\mathrm{exact}}(x)=\frac{x(1-x)}{2k_{\mathrm{hom}}}$. Numerically, the LP solution agrees closely with the parabola, with the pointwise relative error remaining uniformly at the $10^{-5}$ level across the domain.

\begin{figure}[H]
\centering
\includegraphics[width=0.8\textwidth]{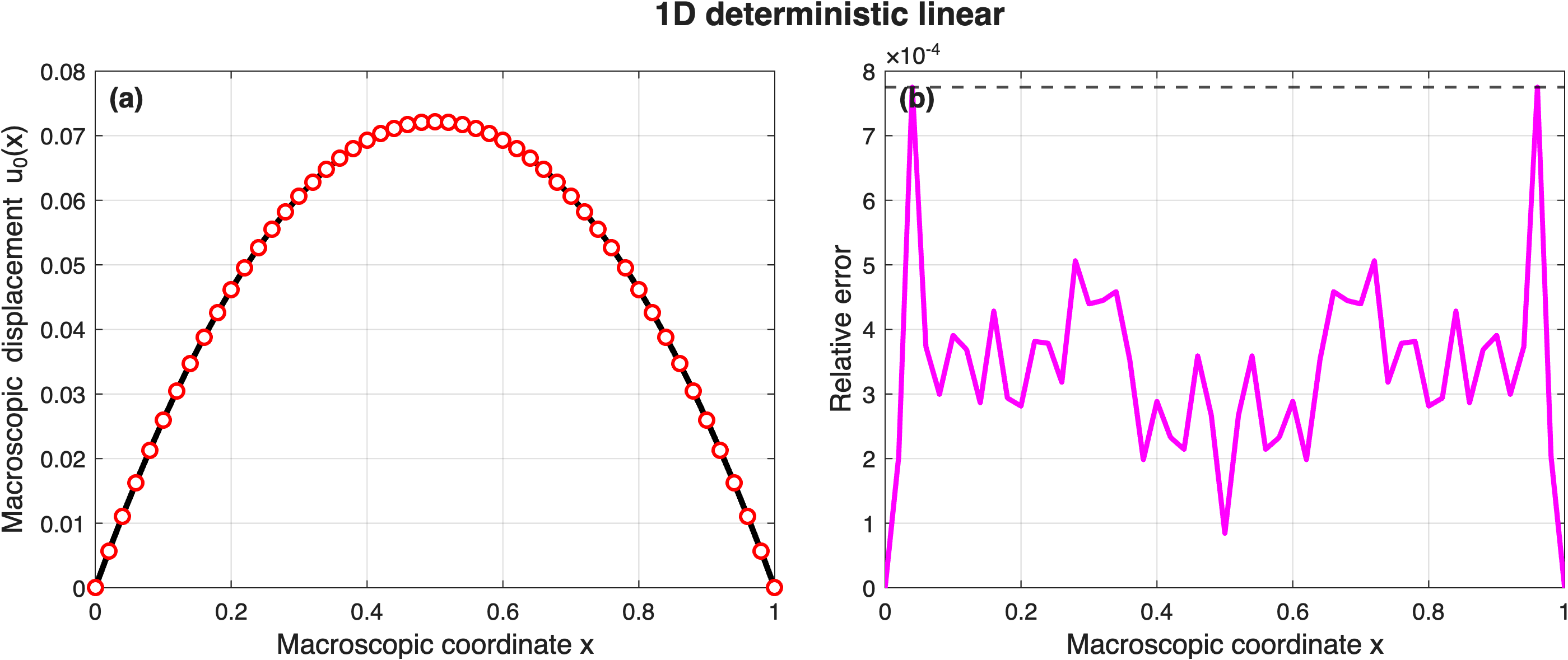}
\caption{1D deterministic linear benchmark. (a): exact solution (black line) and LP numerical solution (red dots). (b): pointwise relative-error profile.}
\end{figure}

\begin{figure}[H]
\centering
\includegraphics[width=0.8\textwidth]{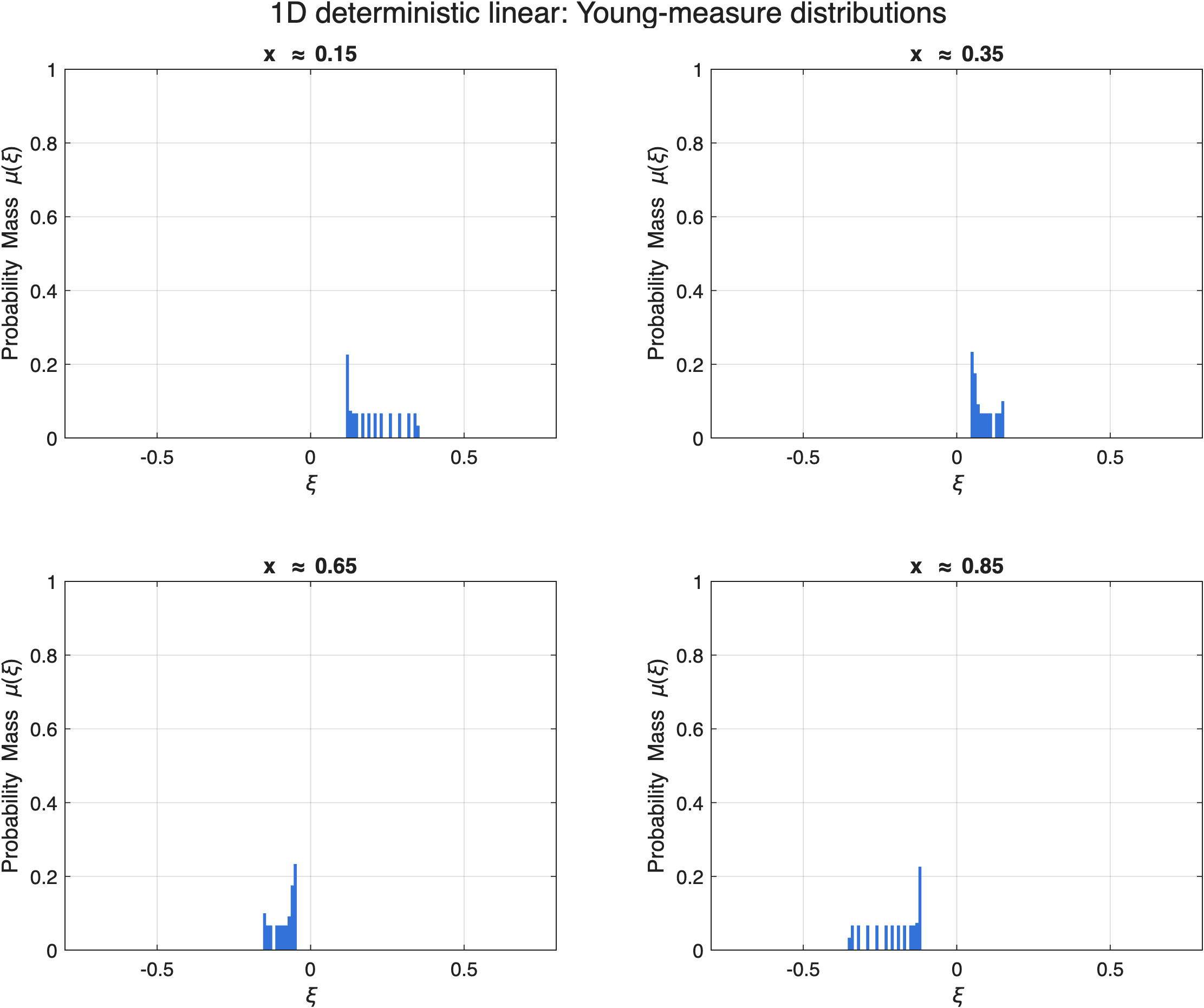}
\caption{1D deterministic linear: marginal Young-measure distribution at representative macroscopic locations (bars show $\int_Y\mu_{x,y}(\xi)\,dy$).}
\end{figure}

\paragraph{Discussion.}
In the linear case, the Young measure is theoretically a Dirac mass $\delta_{\xi^*(x,y)}$ at each $(x,y)$, where $\xi^*(x,y) = u_0'(x) + w_y'(y,u_0'(x))$ is determined by the cell corrector $w$. The marginal Young-measure distribution shown in Figure~2 therefore reflects the pushforward of Lebesgue measure on $Y$ through the analytic map $y\mapsto\xi^*(x,y)$: it is smooth and unimodal, concentrated near the macroscopic gradient $u_0'(x)$ at each $x$. Since $k(y)=2+\sin(2\pi y)$ is analytic, both $\xi^*(x,y)$ and the effective coefficient $k_{\rm hom}=\sqrt{3}$ are determined by analytic corrector computations. The excellent agreement at the $10^{-5}$ level confirms that the LP correctly recovers this analytic structure through the discrete gradient-state representation.

\subsubsection{Nonlinear case}
\begin{table}[H]
\centering
\small
\begin{tabular}{@{}ll@{}}
\toprule
Parameter & Value \\
\midrule
Macro grid & $50$ \\
Micro grid & $30$ \\
State grid $\xi$ & $201$, range $[-1,1]$ \\
Constitutive law & $W=\frac{1}{4}k(y)\xi^4$, $a=k(y)\xi^3$ (cubic flux) \\
\bottomrule
\end{tabular}
\caption{1D deterministic nonlinear: LP discretization setup.}
\label{tab:disc_1d_det_nonlin}
\end{table}

The deterministic nonlinear benchmark is
\begin{equation}
-\frac{d}{dx}\left(k\!\left(\frac{x}{\varepsilon}\right)u_\varepsilon'(x)^3\right)=1,
\qquad u_\varepsilon(0)=u_\varepsilon(1)=0,
\end{equation}
with the same $k(y)=2+\sin(2\pi y)$. The effective coefficient is
\begin{equation}
k_{\mathrm{hom}}=\left(\int_0^1k(y)^{-1/3}\,dy\right)^{-3},
\end{equation}
and the exact homogenized profile is
\begin{equation}
u_{\mathrm{exact}}(x)=\frac{3}{4}k_{\mathrm{hom}}^{-1/3}\left(2^{-4/3}-\left|x-\tfrac12\right|^{4/3}\right).
\end{equation}
The computed solution captures the non-quadratic shape accurately.

\begin{figure}[H]
\centering
\includegraphics[width=0.9\textwidth]{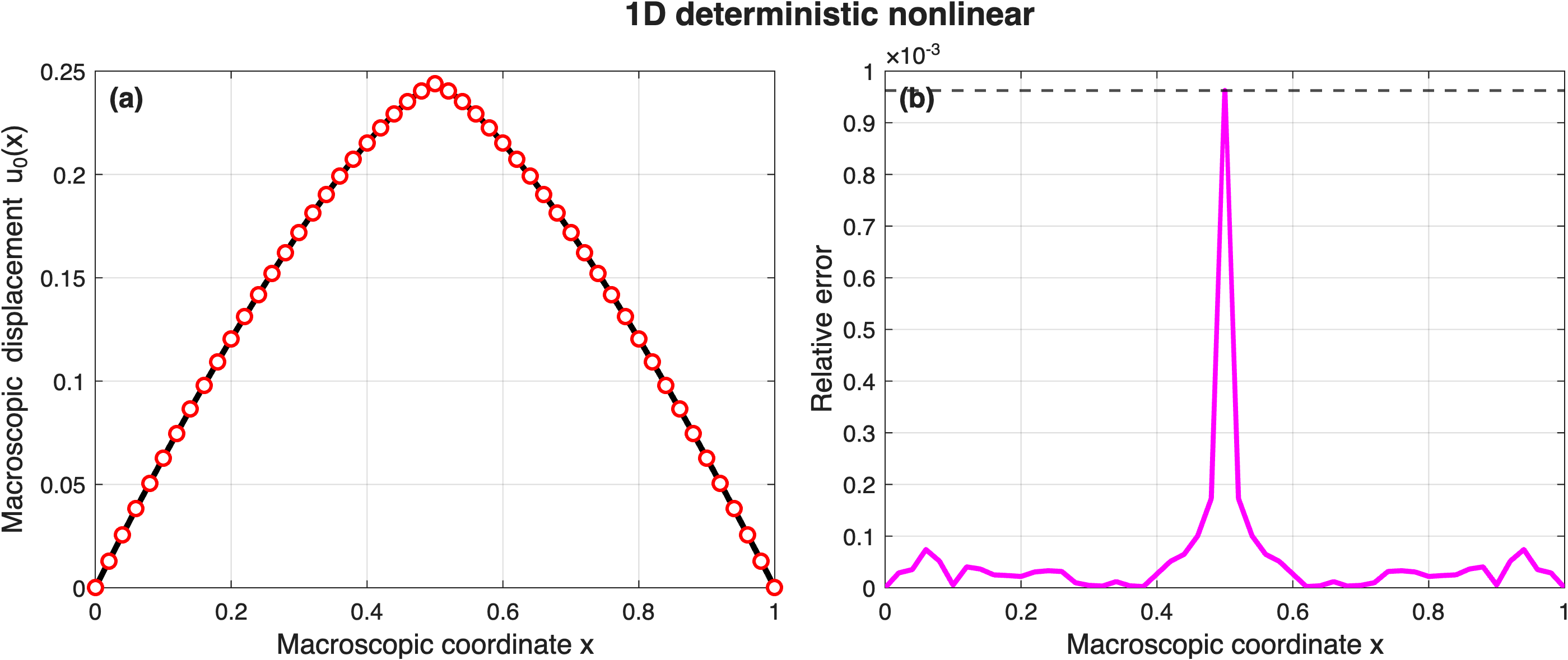}
\caption{1D deterministic nonlinear benchmark. Panel (a): exact solution (black line) and LP numerical solution (red dots). Panel (b): pointwise relative-error profile.}
\end{figure}

\begin{figure}[H]
\centering
\includegraphics[width=0.8\textwidth]{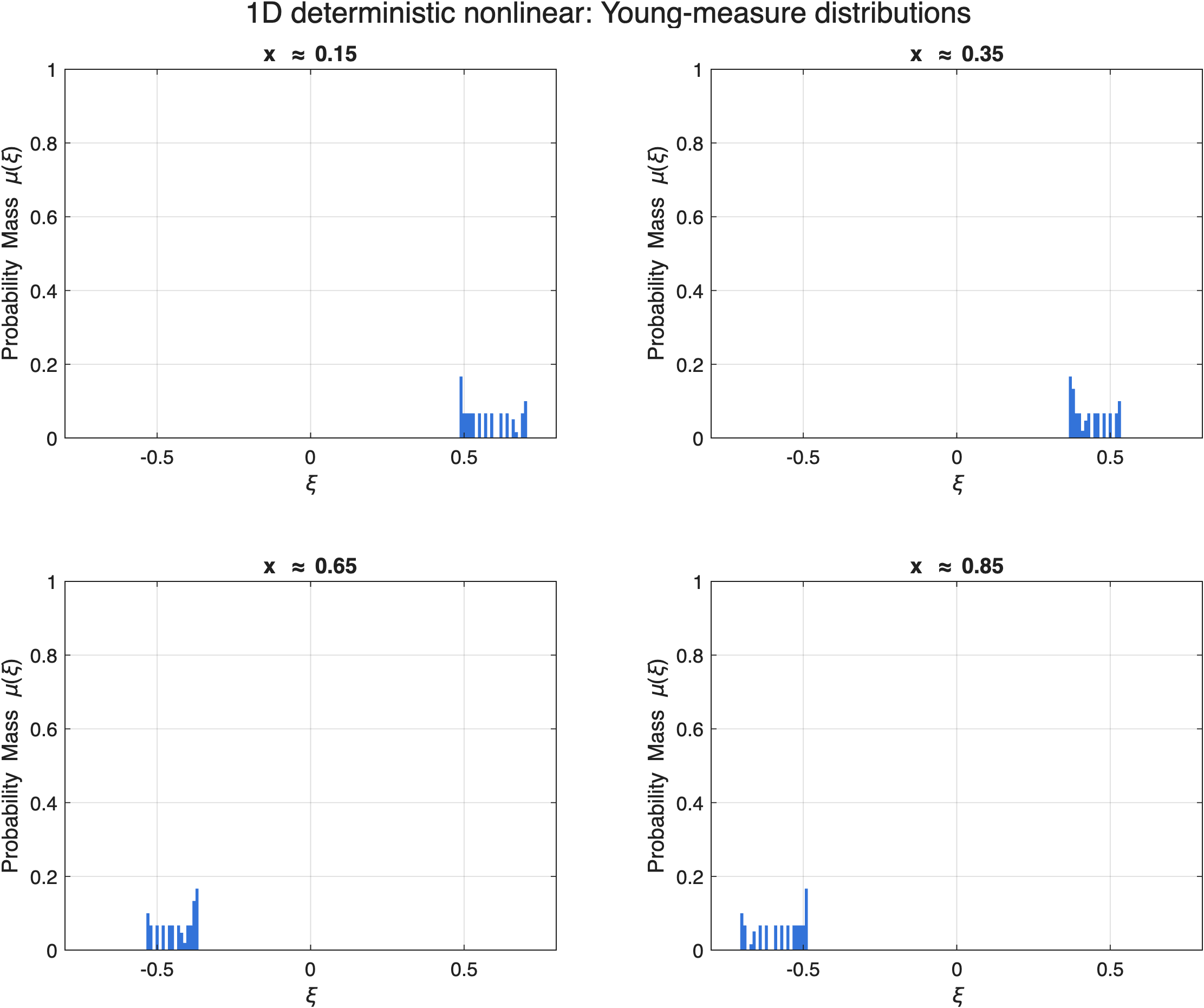}
\caption{1D deterministic nonlinear: marginal Young-measure distribution at representative macroscopic locations (bars show $\int_Y\mu_{x,y}(\xi)\,dy$).}
\end{figure}

\paragraph{Discussion.}
The nonlinear cubic-flux operator $a(y,\xi)=k(y)\xi^3$ satisfies the strict monotonicity condition $\partial_\xi a>0$ for $k>0$, so the cell corrector map $y\mapsto\xi^*(x,y)$ is a diffeomorphism on $Y$ for each fixed macroscopic gradient. Accordingly, the marginal Young measure (Figure~4) is again the pushforward of Lebesgue measure on $Y$ through this map and appears as a smooth, spread distribution in $\xi$. The non-quadratic profile $u_{\rm exact}\propto|x-1/2|^{4/3}$ has a Hölder singularity at $x=1/2$, where $u_{\rm exact}'=0$; this causes a mild degradation in the spatial convergence rate near the center of the domain (visible in the error panel), consistent with the algebraic $O(N_x^{-4/3})$ rate expected for $W^{1,\infty}$-approximation of a $C^{4/3}$ function. The LP nonetheless captures the overall shape with errors below $10^{-3}$, demonstrating robustness to mild solution singularities.

\subsection{2D Deterministic Periodic Media}
For both 2D tests,
\begin{equation}
\Omega=(0,1)^2,
\qquad u=0\text{ on }\partial\Omega,
\end{equation}
with anisotropic periodic coefficients
\begin{equation}
k_1(s)=2+\sin(2\pi s),\qquad k_2(s)=3+\cos(2\pi s),
\end{equation}
and manufactured exact field
\begin{equation}
u_{\mathrm{exact}}(x_1,x_2)=x_1(1-x_1)x_2(1-x_2).
\end{equation}

\subsubsection{Linear anisotropic case}
\begin{table}[H]
\centering
\small
\begin{tabular}{@{}ll@{}}
\toprule
Parameter & Value \\
\midrule
Macro grid & $20\times 20$  \\
Micro grid & $5\times 5$  \\
State grid & $25\times 25$  \\
Micro coefficients & $k_1(y_1)=2+\sin(2\pi y_1)$, $k_2(y_2)=3+\cos(2\pi y_2)$ \\
Constitutive law & $W(y,\xi)=\frac{1}{2}(k_1\xi_1^2+k_2\xi_2^2)$ (linear) \\
\bottomrule
\end{tabular}
\caption{2D linear anisotropic: LP discretization setup.}
\label{tab:disc_2d_lin}
\end{table}

The oscillatory equation is
\begin{equation}
-\nabla\cdot\left(\begin{bmatrix}k_1(x_1/\varepsilon)&0\\0&k_2(x_2/\varepsilon)\end{bmatrix}\nabla u_\varepsilon\right)=f_{\mathrm{lin}}(x),
\qquad u_\varepsilon|_{\partial\Omega}=0.
\end{equation}
The right-hand side is chosen as
\begin{equation}
f_{\mathrm{lin}}(x_1,x_2)=2k_{1,\mathrm{eff}}x_2(1-x_2)+2k_{2,\mathrm{eff}}x_1(1-x_1),
\end{equation}
where $k_{1,\mathrm{eff}}$ and $k_{2,\mathrm{eff}}$ are the homogenized coefficients.

\begin{figure}[H]
\centering
\includegraphics[width=0.95\textwidth]{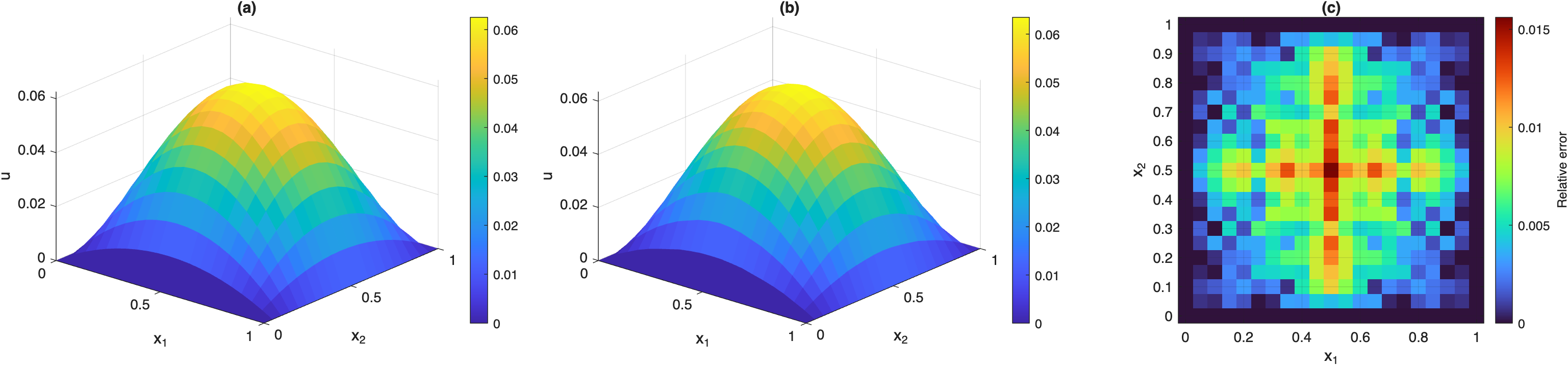}
\caption{2D linear case: (a) exact field, (b) LP-computed field, and (c) relative-error heat map, with a maximum relative error of $1.561\times10^{-2}$.}
\end{figure}

\begin{figure}[H]
\centering
\includegraphics[width=0.62\textwidth]{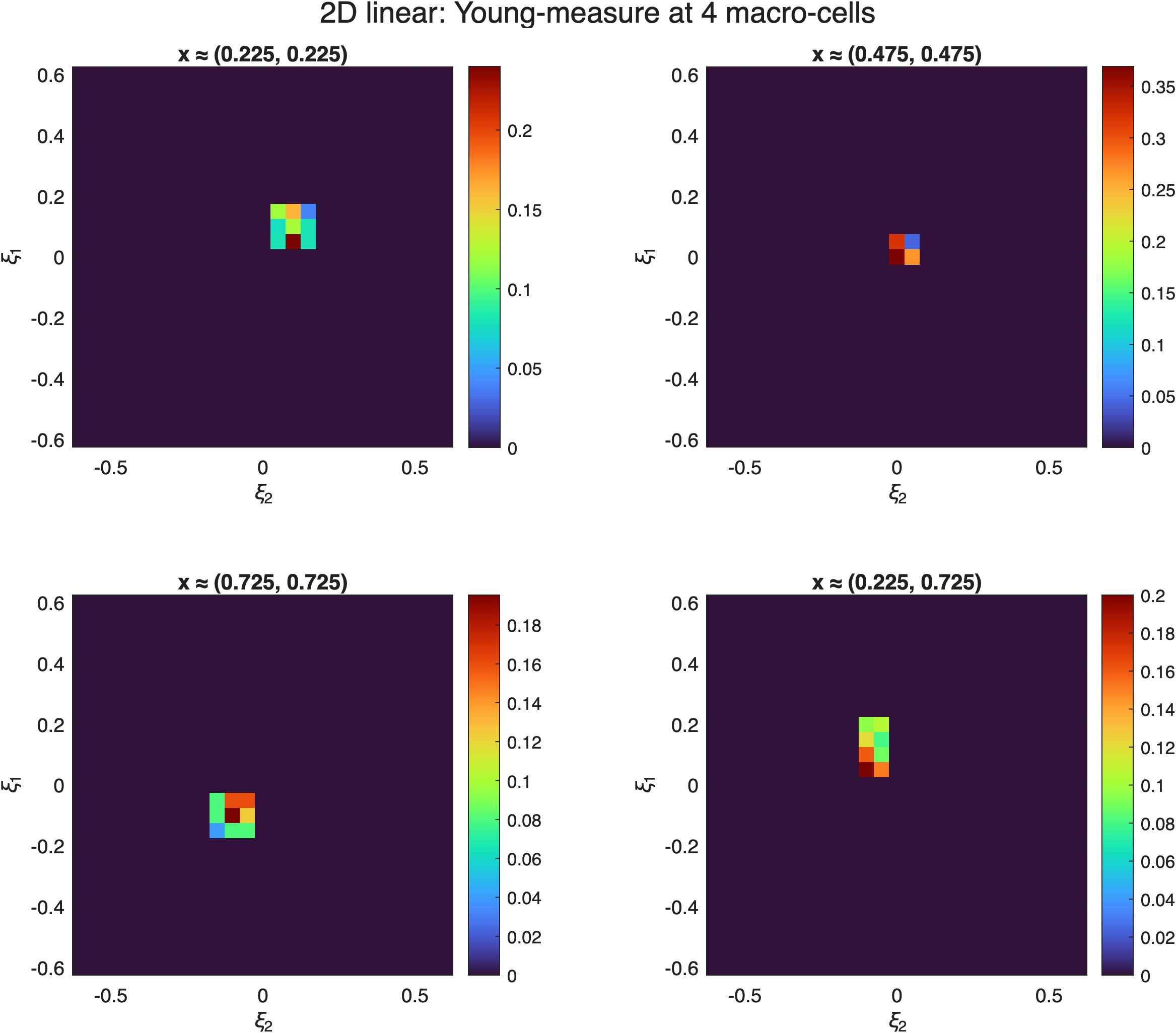}
\caption{2D linear case: marginal Young-measure distribution at 4 representative central macro-cells. The horizontal and vertical axes are $\xi_1$ and $\xi_2$, and the plotted quantity is $\int_Y\mu_{x,y}(\xi_1,\xi_2)\,dy$.}
\end{figure}

\paragraph{Discussion.}
The anisotropic coefficient structure decouples the two spatial directions: the effective coefficients $k_{1,\rm eff}=\sqrt{3}$ and $k_{2,\rm eff}=1/(\int_0^1(3+\cos 2\pi s)^{-1}ds)$ are determined independently by 1D cell problems in $y_1$ and $y_2$ respectively. The marginal Young-measure distributions (Figure~6) reflect this decoupling: the support in the $(\xi_1,\xi_2)$-plane is concentrated along a curve determined by the product of two independent corrector maps. The maximum relative error of $1.561\times10^{-2}$ is consistent with the coarser 2D discretization ($20\times20$ macro grid, $5\times5$ micro grid, $25\times25$ gradient grid), where the LP dimension $n\sim(20)^2(5)^2(25)^2\approx 1.56\times10^7$ is already large for classical solvers. This illustrates the LP dimension scaling $n\sim N_x^d N_y^d N_\xi^d$ in 2D and motivates the use of quantum LP solvers for larger problems.

\subsubsection{Nonlinear anisotropic case}
\begin{table}[H]
\centering
\small
\begin{tabular}{@{}ll@{}}
\toprule
Parameter & Value \\
\midrule
Macro grid & $15\times 15$  \\
Micro grid & $5\times 5$  \\
State grid & $35\times 35$  \\
Micro coefficients & $k_1(y_1)=2+\sin(2\pi y_1)$, $k_2(y_2)=3+\cos(2\pi y_2)$ \\
Constitutive law & $W(y,\xi)=\frac{1}{4}(k_1\xi_1^4+k_2\xi_2^4)$ (cubic flux) \\
\bottomrule
\end{tabular}
\caption{2D nonlinear anisotropic: LP discretization setup.}
\label{tab:disc_2d_nonlin}
\end{table}

The oscillatory equation is
\begin{equation}
-\partial_{x_1}\!\left(k_1(x_1/\varepsilon)(\partial_{x_1}u_\varepsilon)^3\right)
-\partial_{x_2}\!\left(k_2(x_2/\varepsilon)(\partial_{x_2}u_\varepsilon)^3\right)
=f_{\mathrm{nonlin}}(x),
\qquad u_\varepsilon|_{\partial\Omega}=0.
\end{equation}

The right-hand side is chosen as
\begin{equation}
f_{\mathrm{nonlin}}(x_1,x_2)=6K_{1,\mathrm{eff}}x_2^3(1-x_2)^3(1-2x_1)^2+6K_{2,\mathrm{eff}}x_1^3(1-x_1)^3(1-2x_2)^2,
\end{equation}
where $K_{1,\mathrm{eff}}$ and $K_{2,\mathrm{eff}}$ are the homogenized coefficients.

\begin{figure}[H]
\centering
\includegraphics[width=0.92\textwidth]{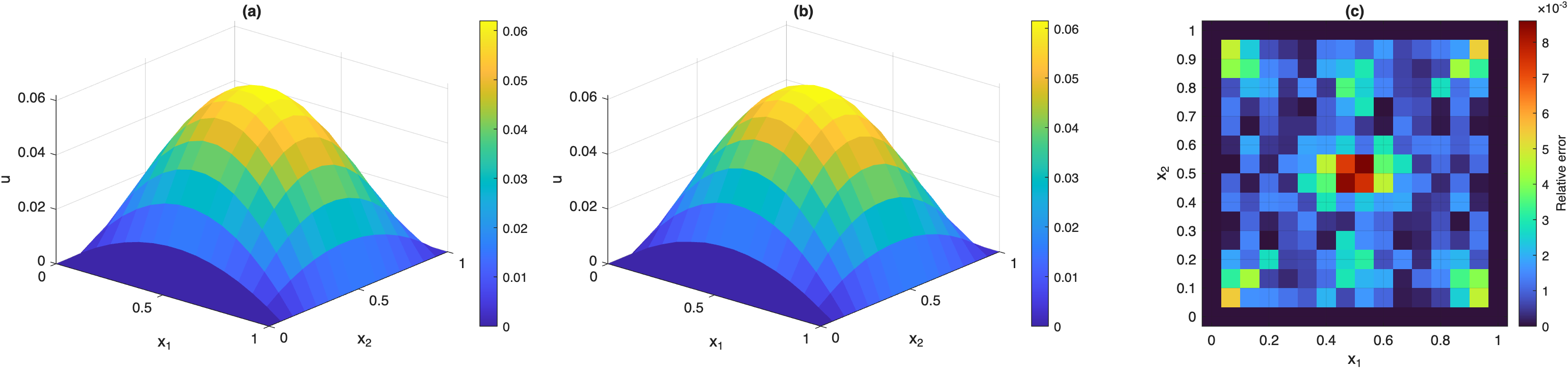}
\caption{2D nonlinear case: (a) exact field, (b) LP-computed field, and (c) relative-error heat map, with a maximum relative error of $8.599\times10^{-3}$.}
\end{figure}

\begin{figure}[H]
\centering
\includegraphics[width=0.6\textwidth]{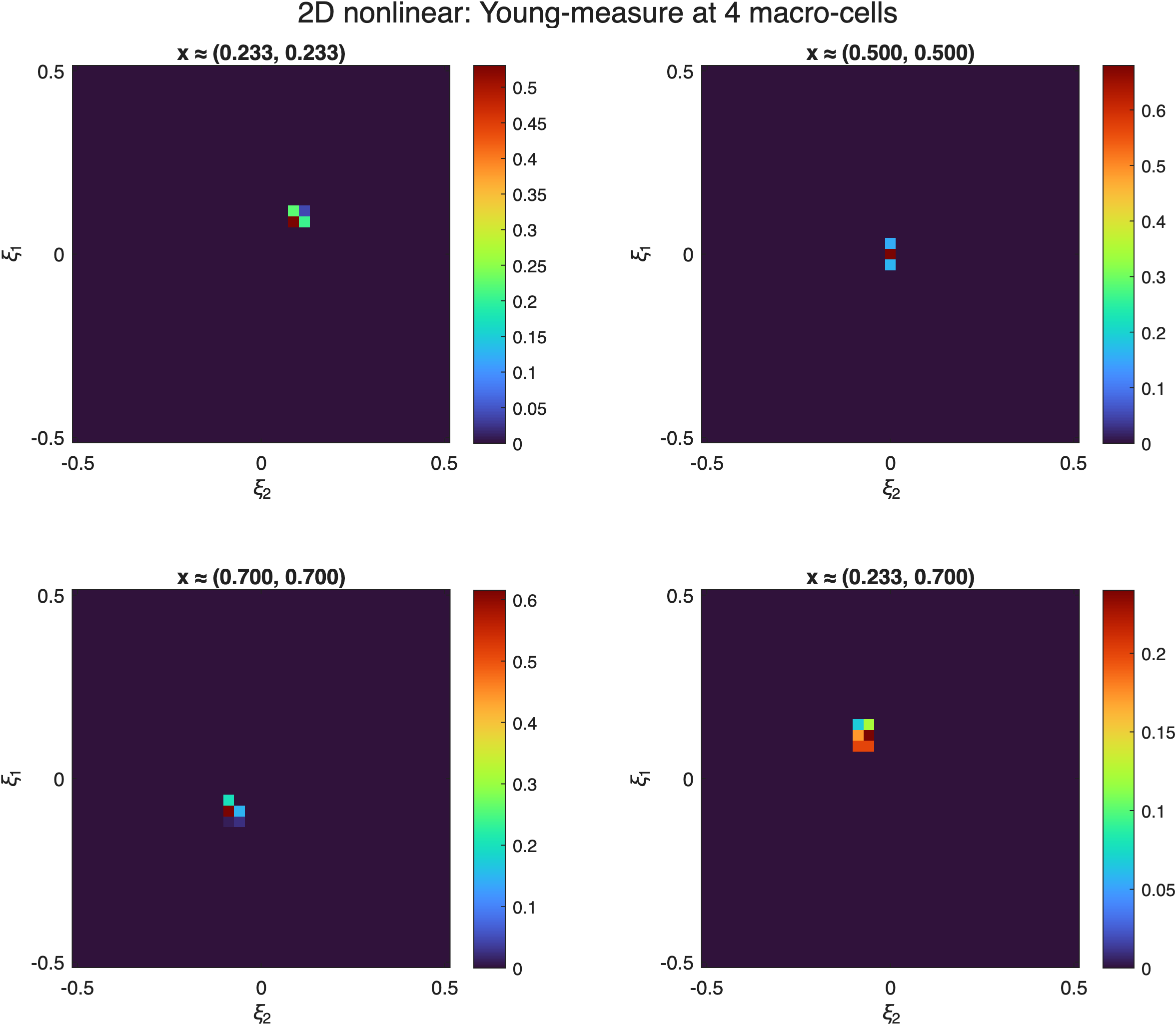}
\caption{2D nonlinear case: marginal Young-measure distribution at 4 representative central macro-cells. The horizontal and vertical axes are $\xi_1$ and $\xi_2$, and the plotted quantity is $\int_Y\mu_{x,y}(\xi_1,\xi_2)\,dy$.}
\end{figure}

\paragraph{Discussion.}
The 2D nonlinear case combines the decoupled anisotropic structure with a cubic-flux nonlinearity in each direction. The effective coefficients $K_{i,\rm eff}=(k_{i,\rm eff}^{-1/3})^{-3}$ are determined by the nonlinear 1D cell problems and appear as the independent cubic homogenization constants for each direction. The marginal Young-measure distributions (Figure~8) are more spread in the $\xi$-plane compared to the linear case, reflecting the stronger nonlinear distortion of the corrector map $\xi^*(x,y)$ under the cubic flux. The maximum relative error of $8.599\times10^{-3}$ is achieved on a $15\times15$ macro grid with $35\times35$ gradient states, and is smaller than the 2D linear error despite the coarser macro grid, because the nonlinear effective profile $u_{\rm exact}=x_1(1-x_1)x_2(1-x_2)$ is smooth and well-resolved by the LP. The result confirms that the Young-measure LP handles nonlinear constitutive laws in multiple dimensions without structural modification.

\subsection{1D Stochastic Media}
\subsubsection{Random linear case}
\begin{table}[H]
\centering
\small
\begin{tabular}{@{}ll@{}}
\toprule
Parameter & Value \\
\midrule
Macro grid & $51$ \\
Random states & $N_\omega=2$: $\{c_1=5, c_2=1\}$ with $\pi=[0.4,0.6]$ \\
State grid $\xi$ & $201$ points, range $[-3,3]$ \\
Constitutive law & linear flux $a(\omega)\xi$ \\
\bottomrule
\end{tabular}
\caption{1D random linear: LP discretization setup.}
\label{tab:disc_1d_rand_lin}
\end{table}

We solve
\begin{equation}
-\frac{d}{dx}\left(c\!\left(\frac{x}{\varepsilon},\omega\right)u_\varepsilon'(x,\omega)\right)=1,
\qquad u_\varepsilon(0,\omega)=u_\varepsilon(1,\omega)=0,
\end{equation}
with $c\in\{5,1\}$ and probabilities $[0.4,0.6]$. The effective coefficient is
\begin{equation}
c_{\mathrm{hom}}=\left(\sum_i p_i/c_i\right)^{-1},
\end{equation}
and $u_{\mathrm{exact}}(x)=x(1-x)/(2c_{\mathrm{hom}})$. Script output gives
\begin{equation}
\max|u_{\mathrm{num}}-u_{\mathrm{exact}}|=7.36\times10^{-4},
\qquad
\frac{\|u_{\mathrm{num}}-u_{\mathrm{exact}}\|_\infty}{\|u_{\mathrm{exact}}\|_\infty}=8.66\times10^{-3}.
\end{equation}

\begin{figure}[H]
\centering
\includegraphics[width=0.92\textwidth]{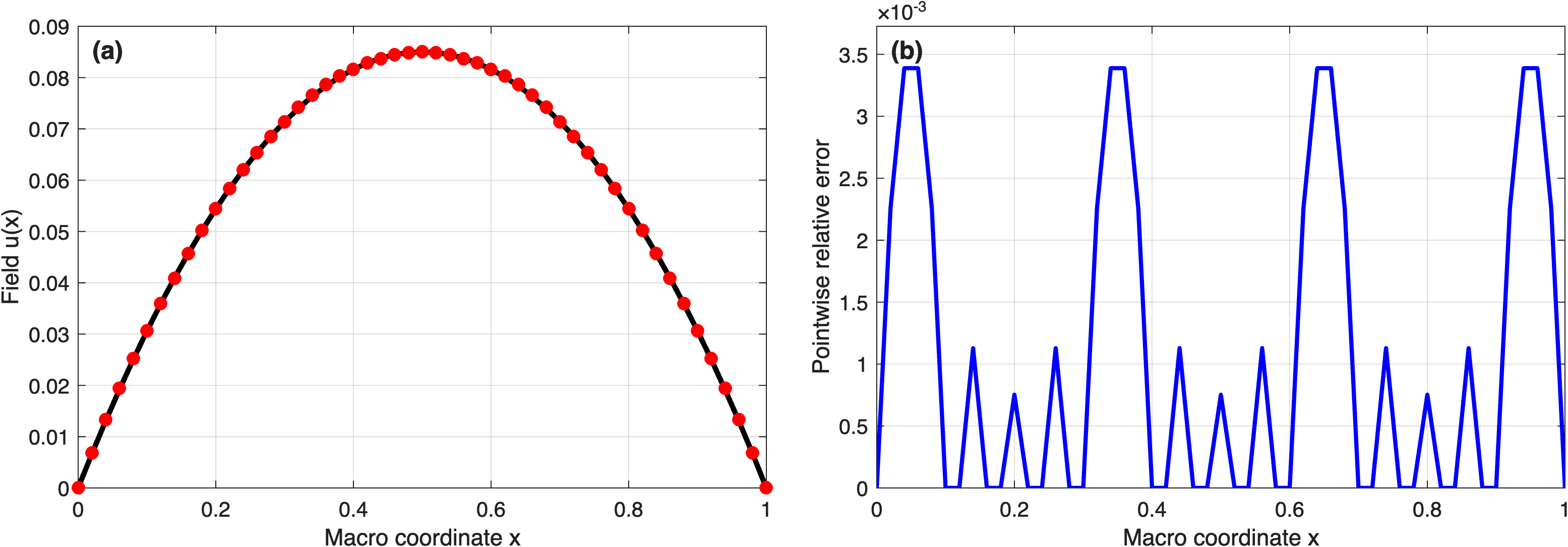}
\caption{1D random linear benchmark. Panel (a): exact solution (black) and LP solution (red). Panel (b): pointwise relative error profile.}
\end{figure}

\begin{figure}[H]
\centering
\includegraphics[width=0.8\textwidth]{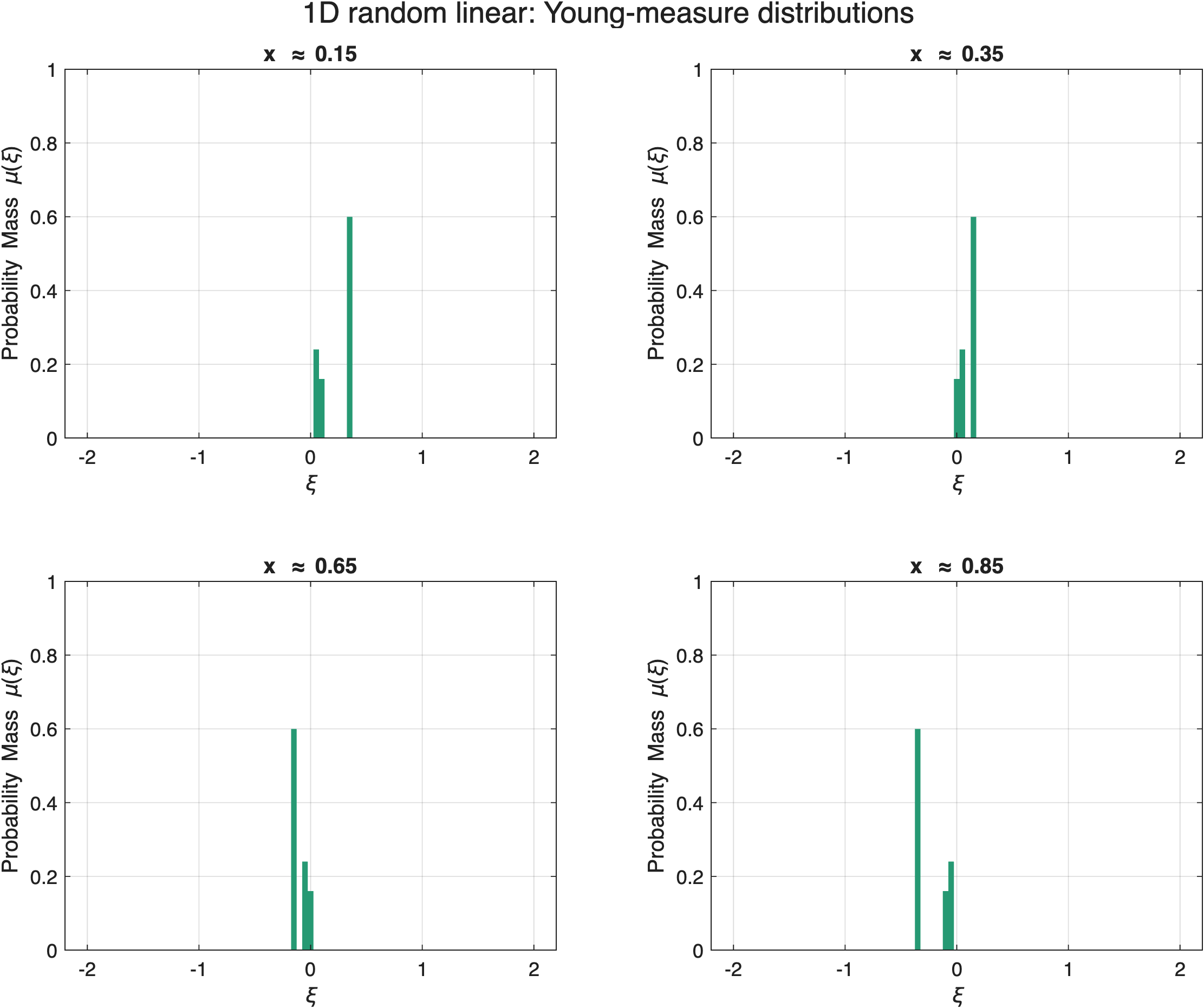}
\caption{1D random linear benchmark: global marginal Young-measure distribution at $4$ representative macroscopic points . The plotted mass is $\hat{\mu}_x(\xi)=\sum_{i=1}^{N_\omega}\mu_{x,i}(\xi)=\sum_{i=1}^{N_\omega}p_i\nu_{x,i}(\xi)$.}
\end{figure}

\paragraph{Discussion.}
With only two random states $\omega\in\{c_1,c_2\}$, the stochastic Young measure $\nu_{x,i}(\xi)$ for each realization $i$ is a Dirac mass $\delta_{\xi^*(x,i)}$ at the corrector gradient for that realization. The global marginal $\hat\mu_x(\xi)=\sum_i p_i\delta_{\xi^*(x,i)}$ is therefore a weighted combination of two Dirac masses, visible as two discrete peaks in Figure~10. The LP correctly identifies both peaks and their weights, recovering the harmonic average effective coefficient $c_{\rm hom}$ to within relative error $8.66\times10^{-3}$. This stochastic case with $N_\omega=2$ represents the simplest instance of the large-$r$ advantage regime: the LP encodes both realizations simultaneously in a single solve, at half the DCS cost of solving two independent fine-scale problems.

\subsubsection{Random nonlinear case}
\begin{table}[H]
\centering
\small
\begin{tabular}{@{}ll@{}}
\toprule
Parameter & Value \\
\midrule
Macro grid  & $51$ \\
Random states & $N_\omega=3$: $\{c_1=1, c_2=3, c_3=6\}$ with $\pi=[0.5,0.3,0.2]$ \\
State grid $\xi$ & $201$ points, range $[-3,3]$ \\
Constitutive law & quadratic flux $a(\omega)|\xi|\xi$ \\
\bottomrule
\end{tabular}
\caption{1D random nonlinear: LP discretization setup.}
\label{tab:disc_1d_rand_nonlin}
\end{table}

We solve
\begin{equation}
-\frac{d}{dx}\left(c\!\left(\frac{x}{\varepsilon},\omega\right)|u_\varepsilon'(x,\omega)|u_\varepsilon'(x,\omega)\right)=1,
\qquad u_\varepsilon(0,\omega)=u_\varepsilon(1,\omega)=0,
\end{equation}
with three random states $c\in\{1,3,6\}$ and probabilities $[0.5,0.3,0.2]$. The effective coefficient is
\begin{equation}
c_{\mathrm{hom}}=\left(\sum_i p_i c_i^{-1/2}\right)^{-2},
\end{equation}
and the exact profile is
$u_{\mathrm{exact}}(x)=\frac{2}{3}c_{\mathrm{hom}}^{-1/2}\left(2^{-3/2}-\left|x-\tfrac12\right|^{3/2}\right)$.

Script output gives
\begin{equation}
\max|u_{\mathrm{num}}-u_{\mathrm{exact}}|=1.27\times10^{-4},
\qquad
\frac{\|u_{\mathrm{num}}-u_{\mathrm{exact}}\|_\infty}{\|u_{\mathrm{exact}}\|_\infty}=7.11\times10^{-4}.
\end{equation}

\begin{figure}[H]
\centering
\includegraphics[width=0.85\textwidth]{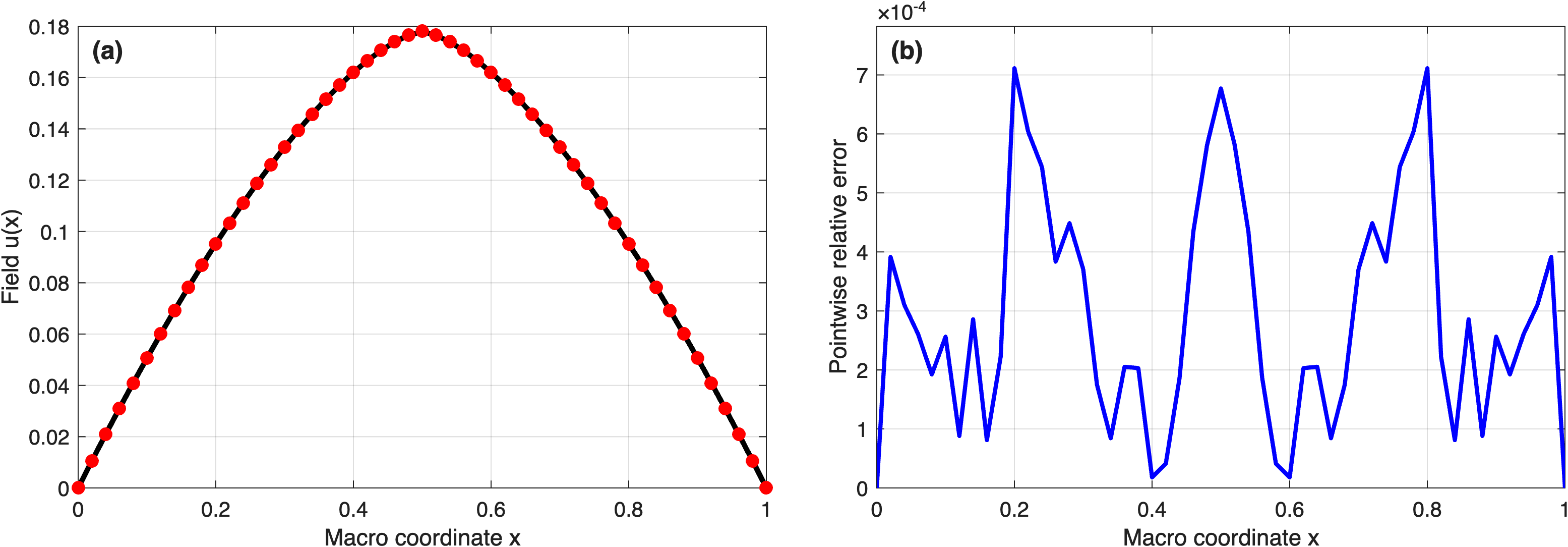}
\caption{1D random nonlinear benchmark. Panel (a): exact solution (black) and LP solution (red). Panel (b): pointwise relative error profile.}
\end{figure}

\begin{figure}[H]
\centering
\includegraphics[width=0.8\textwidth]{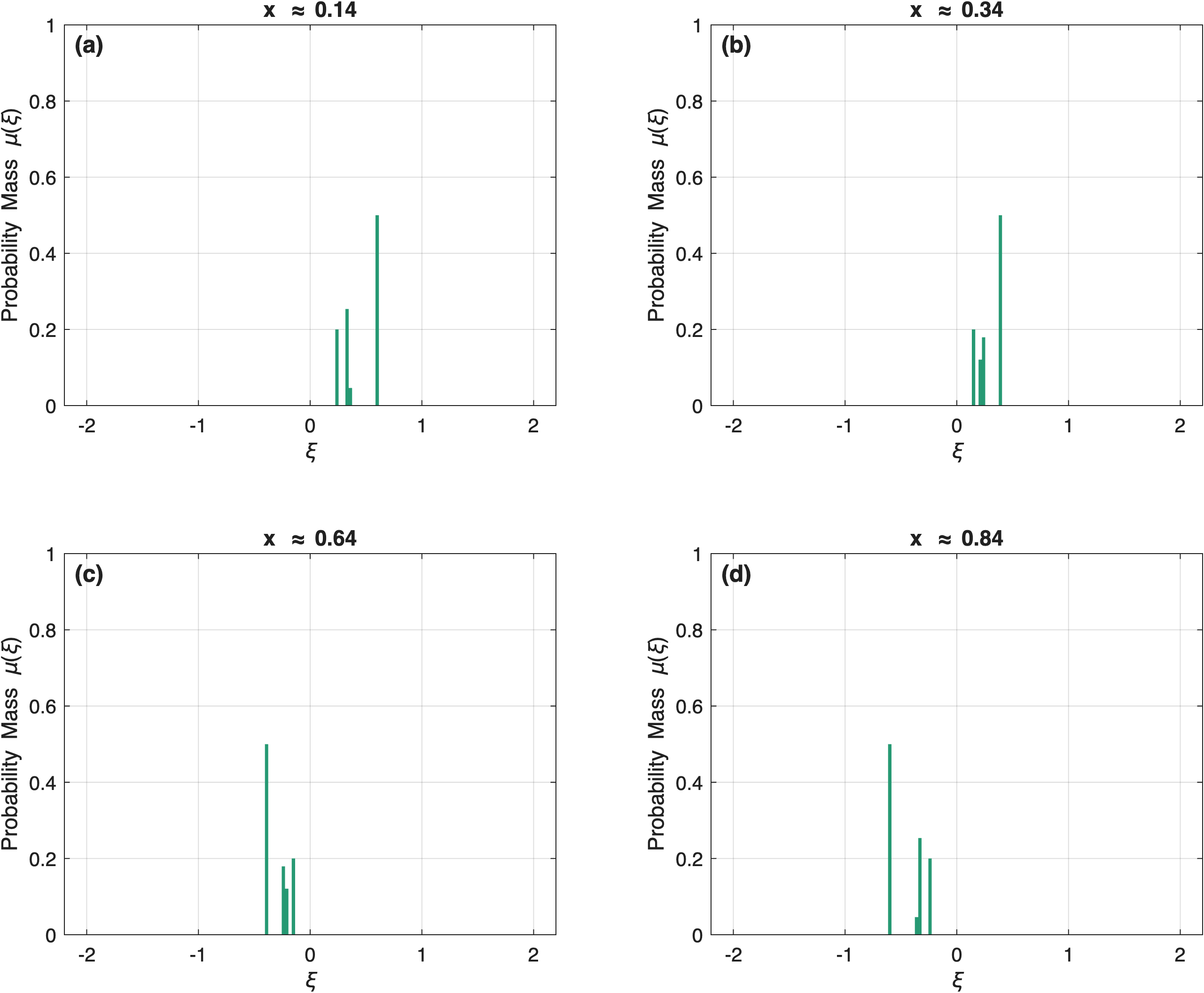}
\caption{1D random nonlinear benchmark: global marginal Young-measure distribution at $4$ representative macroscopic points. In this 1D random implementation (without explicit $y$), the plotted density is $\hat{\mu}_x(\xi)=\sum_{i=1}^{N_\omega}\mu_{x,i}(\xi)=\sum_{i=1}^{N_\omega}p_i\nu_{x,i}(\xi)$.}
\end{figure}

\paragraph{Discussion.}
With three random states and a quadratic flux $a(\omega)|\xi|\xi$, the LP must simultaneously represent three distinct corrector gradients weighted by their probabilities $[0.5,0.3,0.2]$. The marginal Young measure (Figure~12) shows three discrete peaks at $\xi$-values determined by $\xi^*(x,i)\propto c_i^{-1/2}\,{\rm sgn}(u_0'(x))$, each with height proportional to $p_i$. The effective coefficient $c_{\rm hom}=(\sum_i p_i c_i^{-1/2})^{-2}$ blends the three realizations nonlinearly via the cubic homogenization formula. The very small relative error $7.11\times10^{-4}$ demonstrates that the LP correctly resolves the nonlinear stochastic averaging, recovering a relative error nearly an order of magnitude smaller than the linear stochastic case. This improvement is consistent with the smooth $|x-1/2|^{3/2}$ profile being slightly less singular than the linear parabola under the cubic flux energy, and with the finer effective resolution provided by three gradient peaks compared to two.

\subsubsection{Linear non-variational anisotropic case}
\begin{table}[H]
\centering
\small
\begin{tabular}{@{}ll@{}}
\toprule
Parameter & Value \\
\midrule
Macro grid & $14\times 14$  \\
Micro grid & $8\times 8$  \\
State grid & $21\times 21$  \\
Micro coefficients & $a_{11}(y_1)=a_{22}(y_1)=2+\sin(2\pi y_1)$, $a_{12}=1$, $a_{21}=-1$ \\
Constitutive law & $A(y)\xi = \begin{bmatrix} 2+\sin(2\pi y_1) & 1 \\ -1 & 2+\sin(2\pi y_1) \end{bmatrix} \begin{bmatrix} \xi_1 \\ \xi_2 \end{bmatrix}$ (non-variational) \\
\bottomrule
\end{tabular}
\caption{2D linear non-variational anisotropic: LP discretization setup.}
\label{tab:disc_2d_nonvar}
\end{table}

The oscillatory equation is given by
\begin{equation}
-\nabla\cdot\left(A(x/\varepsilon)\nabla u_\varepsilon\right)=f_{\mathrm{nonvar}}(x),
\qquad u_\varepsilon|_{\partial\Omega}=0,
\end{equation}
where the micro-scale coefficient matrix $A(y)$ is explicitly asymmetric. Due to the presence of the skew-symmetric terms ($a_{12}$ and $a_{21}$), the system lacks an underlying scalar energy potential $W(y,\xi)$ such that $\nabla_\xi W = A(y)\xi$. Consequently, the problem is strictly non-variational.

To uniquely determine the homogenized solution within the Young-measure LP framework, the objective function is set to zero (a pure feasibility problem), and the micro-scale divergence-free condition must be explicitly enforced as an equality constraint:
\begin{equation}
    -\nabla_y \cdot \left( A(y)\xi \right) = 0,
\end{equation}
 along with the microscopic curl-free constraint 
 \begin{equation}
     \nabla_y \times \mathbb{E}[\xi] = 0.
 \end{equation}
 To avoid unnatural non-vanishing measure on the phase-space boundary, we impose an artificial energy regularization objective ($c(\xi) = |\xi|^2$) to enforce the collapse of the Young measure into a deterministic Dirac mass. Given the severe degeneracy caused by the measure collapse on the refined phase-space grid, the global LP system is efficiently solved using the interior-point algorithm of the linprog function in MATLAB.
For the method of manufactured solutions, we design the same exact macroscopic solution as $u_0(x) = x_1(1-x_1)x_2(1-x_2)$. The corresponding right-hand side is chosen as

\begin{equation}
f_{\mathrm{nonvar}}(x_1,x_2)=2\sqrt{3}x_2(1-x_2)+4x_1(1-x_1),
\end{equation}

based on the analytically derived homogenized matrix $A_{\mathrm{eff}} = \begin{bmatrix} \sqrt{3} & 1 \\ -1 & 2 \end{bmatrix}$. 

\begin{figure}[H]
\centering
\includegraphics[width=0.95\textwidth]{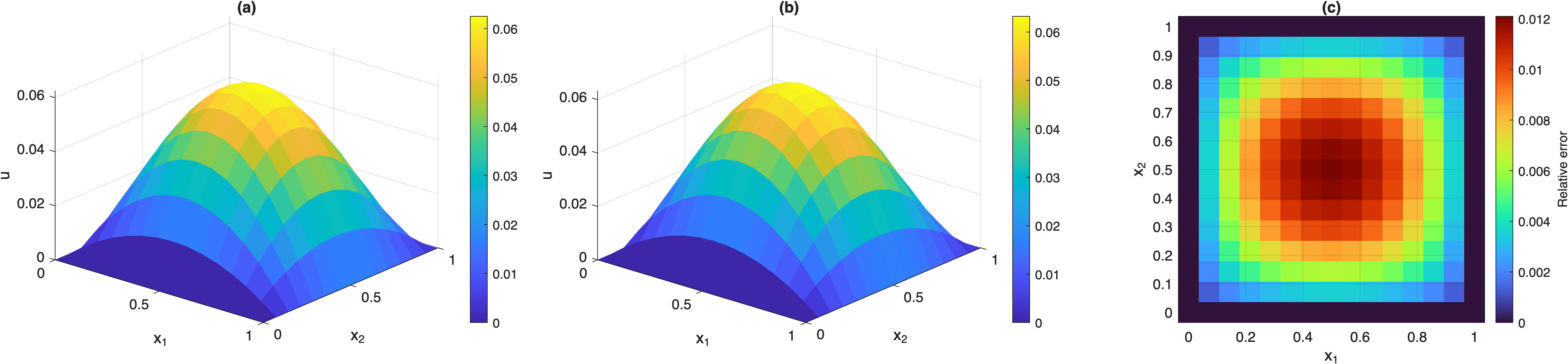}
\caption{Non-variational case: (a) exact field, (b) LP-computed field, and (c) relative-error heat map, with a maximum relative error of $1.2\times10^{-2}$.}
\end{figure}

\begin{figure}[H]
\centering
\includegraphics[width=0.62\textwidth]{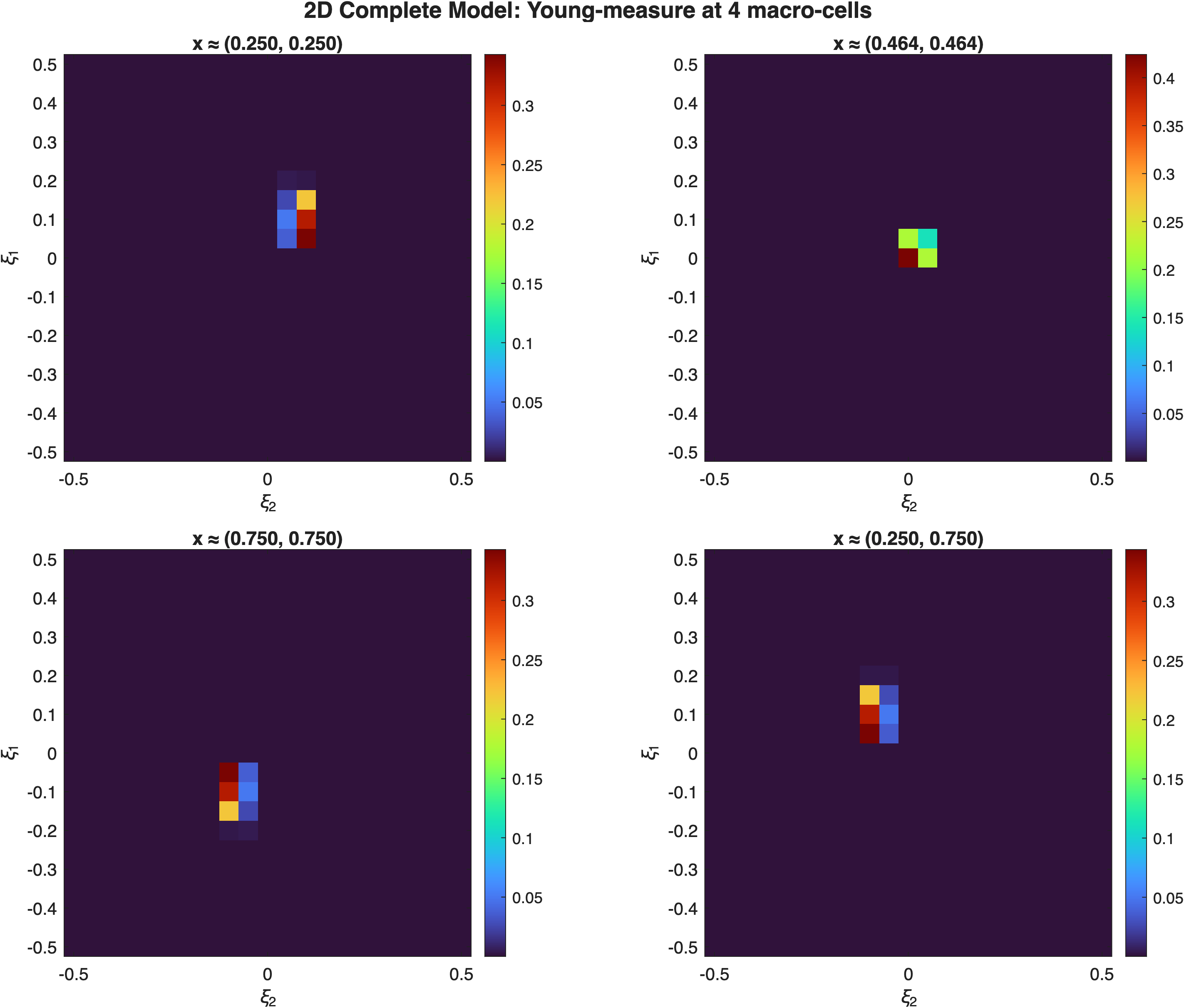}
\caption{Non-variational case: marginal Young-measure distribution at 4 representative central macro-cells. The horizontal and vertical axes are $\xi_1$ and $\xi_2$, and the plotted quantity is $\int_Y\mu_{x,y}(\xi_1,\xi_2)\,dy$.}
\end{figure}

\paragraph{Discussion.}
The non-variational case is qualitatively different from the preceding examples: there is no underlying energy $W(y,\xi)$ whose gradient equals $A(y)\xi$, so the LP cannot be posed as an energy minimization. The feasibility-only formulation with zero objective identifies the unique homogenized solution via the combined microscopic divergence-free and macroscopic curl-free constraints. The energy regularization objective $c(\xi)=|\xi|^2$ is introduced artificially to collapse the Young measure to a Dirac mass, preventing the LP from exploiting degenerate feasible directions in the gradient-state space. The marginal distributions (Figure~14) accordingly concentrate near a single point in $(\xi_1,\xi_2)$ at each macro-cell. The maximum relative error of $1.2\times10^{-2}$ confirms that even in the non-variational setting, the LP feasibility approach correctly reproduces the analytically derived homogenized matrix $A_{\rm eff}=\begin{bmatrix}\sqrt{3}&1\\-1&2\end{bmatrix}$, demonstrating the scope of the Young-measure LP framework beyond energy minimization.

\subsection{Summary on Numerical Results}
The numerical experiments show that the proposed Young-measure LP formulation accurately reproduces known homogenized solutions across deterministic periodic, stochastic, linear, nonlinear, and higher-dimensional elliptic homogenization problems. In the 1D deterministic and stochastic cases, the LP solutions closely match analytical homogenized profiles for both linear and cubic-flux models with very small relative errors. The 2D anisotropic periodic benchmarks further demonstrate agreement with manufactured exact solutions in both linear and nonlinear settings, yielding maximum relative errors of the order of $10^{-2}$, which is quite acceptable given the relatively rough grids limited by the huge number of variables introduced by higher dimensions. Moreover, the non-variational case suggests that our framework also works for problems that do not admit an energy minimization formulation. Overall, these results confirm the correctness, consistency, and flexibility of the Young-measure LP framework for homogenization problems involving periodic media, random media, and nonlinearity.

\section{Conclusion and Outlook}

This work develops a Young-measure linear programming framework for nonlinear deterministic and stochastic homogenization and investigates its compatibility with quantum optimization algorithms. The Young-measure formulation lifts the nonlinear homogenization problem to a linear one in higher dimensions by treating the microscale $y$, the gradient $\xi$, and the random environment $\omega$ as independent variables, thereby characterizing effective macroscopic quantities through a structured LP without directly resolving fine-scale oscillations.

The complexity analysis compares the quantum central path (QCP) algorithm against direct classical solvers (DCS), whose cost $\widetilde O(\varepsilon^{-d})$ is fixed by the need to resolve all $\varepsilon$-scale oscillations regardless of the desired output accuracy. Using first-order ($p=1$) algebraic discretization, two regimes of quantum advantage are identified:
\begin{enumerate}
\item \emph{Deterministic nonlinear homogenization.} Recovering the homogenized solution to accuracy $\delta=\varepsilon^\alpha$ yields QCP cost $\widetilde O(R_1\varepsilon^{-\alpha(3d+4)/2})$, improving over DCS cost $\widetilde O(\varepsilon^{-d})$ whenever $\alpha<2d/(3d+4)$. For $d=2$ this requires $\alpha<2/5$; for $d=3$, $\alpha<6/13$.
\item \emph{Stochastic homogenization at fine-scale accuracy.} With $N_\omega^r$ random discretization points, QCP achieves cost $\widetilde O(R_1\,N_\omega^{r/2}\,\varepsilon^{-(3d/2+1)})$, giving a square-root reduction in stochastic sampling cost that grows with $N_\omega^r$, whenever $N_\omega^r\gtrsim\varepsilon^{-(d+2)}$.
\end{enumerate}
We further conjecture that if the Young-measure density is analytic in $x$ and $y$, spectral discretization in those variables reduces the LP size to $n_{\rm det}\sim\varepsilon^{-d}|\log\varepsilon|^{2d}$ at fine-scale accuracy $\delta=\varepsilon$, yielding QCP cost $\widetilde O(R_1\varepsilon^{-(1+d/2)}|\log\varepsilon|^d)$ and quantum advantage for all $d\ge2$. In the stochastic case, the advantage condition under spectral discretization relaxes to $N_\omega^r\gtrsim\varepsilon^{-(d-2)}$, trivially satisfied for $d\ge2$.

Numerical experiments on 1D and 2D deterministic and stochastic benchmarks confirm the correctness and consistency of the Young-measure LP formulation.

Future work will focus on: rigorous convergence analysis for the Young-measure LP in appropriate function spaces and norms; verification of the spectral conjecture, including analysis of the central-path scale $R_1$ under spectral discretization and of the analytic regularity of the Young-measure density; extensions to broader classes of nonlinear and stochastic PDEs; and comprehensive quantum complexity estimates including state preparation, oracle construction, and measurement costs beyond the query-complexity level. Overall, the proposed framework establishes a bridge among nonlinear homogenization, stochastic multiscale PDEs, linear programming, and quantum scientific computing, and points a promising direction for efficient quantum computation of nonlineaer multiscale and stochastic PDEs.

\section*{Acknowledgments}
SJ and ZL were supported by NSFC grant No. 12341104,  the Shanghai Pilot Program for Basic Research, the Shanghai Jiao Tong University 2030 Initiative,  the  Science and Technology Innovation Key R\&D Program of Chongqing grant No. CSTB2024TIAD-STX0035, and the Fundamental Research Funds for the Central Universities.

\bibliographystyle{plain}
\bibliography{ref}

\end{document}